\documentclass[letterpaper]{article}

\usepackage[margin=2.5cm]{geometry}

\usepackage[leqno]{amsmath}
\usepackage{amsfonts,amssymb,amsthm,amscd,amsxtra}

\usepackage{textcomp}

\usepackage{eucal} 

\usepackage{mathrsfs} 

\let\goth\mathfrak

\usepackage{xspace}

\newtheorem{theorem}{Theorem}[section]
\newtheorem{proposition}[theorem]{Proposition}
\newtheorem{lemma}[theorem]{Lemma}
\newtheorem{corollary}[theorem]{Corollary}

\theoremstyle{definition}
\newtheorem{definition}[theorem]{Definition}
\newtheorem{example}[theorem]{Example}

\theoremstyle{remark}
\newtheorem{remark}[theorem]{Remark}

\numberwithin{equation}{section}

\newcommand{\abs}[1]{\left\lvert#1\right\rvert}

\newcommand{\cuphat}{\,\hat\cup\,}
\newcommand{\dual}[2]{\langle#1,#2\rangle}
\newcommand{\hoeq}[1][]{\underset{#1}{\simeq}} 
\newcommand{\norm}[1]{\lVert#1\rVert}
\newcommand{\tame}[2]{\bigl(#1,#2\bigr]}

\newcommand{\tate}{2\pi\smash{\sqrt{-1}}}
\newcommand{\lAngle}{\langle\!\langle}
\newcommand{\rAngle}{\rangle\!\rangle}
\newcommand{\lto}{\longrightarrow} 
\newcommand{\iso}{\cong} 
\newcommand{\coin}{\equiv} 
\newcommand{\unit}{\times} 

\newcommand{\onehalf}{\frac{1}{2}}
\newcommand{\onefourth}{\frac{1}{4}}
\newcommand{\ihalf}{\frac{\sqrt{-1}}{2}}
\newcommand{\eqdef}{\overset{\text{\normalfont\tiny def}}{=}}
\newcommand{\cm}[1]{\mathscr{CM}(#1)} 
\newcommand{\hyp}{\mathsf{hyp}}
\newcommand{\good}{\mathit{good}}
\newcommand{\pre}{\mathit{pre}}
\newcommand{\bei}{Be\u\i{}linson\xspace}
\newcommand{\cech}{\v{C}ech\xspace}

\newcommand{\field}[1]{\ensuremath{\mathbf{#1}}}
\newcommand{\ZZ}{\field{Z}}

\newcommand{\RR}{\field{R}}
\newcommand{\CC}{\field{C}}

\newcommand{\NN}{\field{N}}

\newcommand{\sheaf}[1]{\mathscr{#1}}

\newcommand{\qi}{\overset{\simeq}{\rightarrow}}
\newcommand{\lqi}{\overset{\simeq}{\longrightarrow}}

\newcommand{\coho}[1]{\mathrm{#1}}%
\renewcommand{\H}{\coho{H}} 

\newcommand{\hypercoho}[1]{\mathbb{#1}}%
\newcommand{\HHH}{\hypercoho{H}}

\newcommand{\cover}[1]{\mathfrak{#1}}
\newcommand{\vC}[3][\bullet]{\Check{\mathrm{C}}^{#1}(#2,#3)}

\newcommand{\complex}[2][\bullet]{{#2}^{#1}}

\newcommand{\CCC}[1][\bullet]{\complex[#1]{\mathsf{C}}}

\DeclareMathOperator{\coker}{Coker} \DeclareMathOperator{\cone}{Cone}
\DeclareMathOperator{\dd}{d}\renewcommand{\d}{\dd\mspace{-2mu}} 
\newcommand{\dc}{\dd^c\mspace{-2mu}} 

\newcommand{\del}{\partial} 
\newcommand{\delb}{\Bar\partial} 
\newcommand{\deltacheck}{\delta} 

\DeclareMathOperator{\Pic}{Pic}
\newcommand{\Pichat}{\widehat\Pic}

\DeclareMathOperator{\supp}{supp}
\DeclareMathOperator{\Tot}{Tot}

\newcommand{\sha}[2][\bullet]{\sheaf{A}_{#2}^{#1}}
\newcommand{\she}[2][\bullet]{\sheaf{E}_{#2}^{#1}}
\newcommand{\sho}[1]{\sheaf{O}_{#1}}
\newcommand{\shomega}[2][\bullet]{\varOmega_{#2}^{#1}}
\newcommand{\deligne}[4][\bullet]{%
  \def\tempa{}%
  \def\tempb{#2}%
  \ifx\tempa\tempb
  #3(#4)^{#1}_\mathcal{D}
  \else
  #3(#4)^{#1}_{\mathcal{D},#2}
  \fi
}
\newcommand{\deltilde}[4][\bullet]{%
  \smash[t]{\widetilde{#3_{#2}(#4)}}^{#1}_{\mathcal{D}}}
\newcommand{\delH}[4][\bullet]{{\H}^{#1}_\mathcal{D}(#2, #3(#4))}
\newcommand{\delz}[2][\bullet]{\deligne[#1]{}{\ZZ}{#2}} 
\newcommand{\delr}[2][\bullet]{\deligne[#1]{}{\RR}{#2}} 

\newcommand{\delalg}[3][\bullet]{\goth{D}^{#1}(#2,#3)}

\newcommand{\subdelalg}[4][\bullet]{\goth{D}^{#1}_{#4}(#2,#3)}

\newcommand{\brydhh}[3][\bullet]{C(#3)^{#1}_{#2}}
\newcommand{\dhh}[3][\bullet]{D_\mathit{h.h.}(#3)^{#1}_{#2}}
\newcommand{\ndhh}[3][\bullet]{
  \def\tempa{}%
  \def\tempb{#2}%
  \ifx\tempa\tempb
  \goth{D}_\mathit{h.h.}(#3)^{#1}
  \else
  \goth{D}_\mathit{h.h.}(#3)^{#1}_{#2}
  \fi
}
\newcommand{\dhhH}[3][\bullet]{{\widehat\H}^{#1}_{\mathcal{D}}(#2;#3)}

\title{Hermitian-holomorphic Deligne cohomology, Deligne pairing
  for singular metrics, and hyperbolic metrics}

\author{Ettore Aldrovandi\\
  Department of Mathematics\\
  Florida State University\\
  Tallahassee, FL 32306-4510, USA\\
  \texttt{aldrovandi@math.fsu.edu} }

\date{}

\begin{document}

\maketitle

\noindent
\thanks{Dedicated to \textmu-T}

\begin{abstract}
  We introduce a model for Hermitian holormorphic Deligne
  cohomology on a projective algebraic manifold which allows to
  incorporate singular hermitian structures along a normal
  crossing divisor. In the case of a projective curve, the
  cup-product in cohomology is shown to correspond to a
  generalization of the Deligne pairing to line bundles with
  ``good'' hermitian metrics in the sense of Mumford and others.
  A particular case is that of the tangent bundle of the curve
  twisted by the negative of the singularity divisor of a
  hyperbolic metric: its cup square (corrected by the total area)
  is shown to be a functional whose extrema are the metrics of
  constant negative curvature.
\end{abstract}

\tableofcontents

\section{Introduction}
\label{sec:introduction}

\subsection{Preliminaries and statement of the results}
\label{sec:prel-stat-results}

Let $X$ be a smooth projective curve over $\CC$ (a Riemann
surface) and let $\sheaf{L}$ and $\sheaf{M}$ be two line bundles
($=$ invertible sheaves) over $X$. In \cite{MR89b:32038} Deligne
defines a pairing $\dual{\sheaf{L}}{\sheaf{M}}$ as a complex line
(i.e. a one-dimensional vector space over \CC). Moreover, if both
line bundles carry smooth hermitian fiber metrics $\norm{\cdot}$,
then the complex line carries a hermitian scalar product, whose
expression was also defined in \cite{MR89b:32038} in terms of two
conveniently chosen rational sections $l$ and $m$. (See
section~\ref{sec:determ-cohom} for the relevant formulas.)

A remarkable fact is that the metrized line
$(\dual{\sheaf{L}}{\sheaf{M}}, \norm{\cdot})$ can be obtained by
entirely cohomological means, as a cup product in Hermitian
holomorphic Deligne cohomology between the classes of the
metrized line bundles $\Bar{\sheaf{L}} = (\sheaf{L},
\norm{\cdot})$ and $\Bar{\sheaf{M}} = (\sheaf{M}, \norm{\cdot})$,
\cite{math.CV/0211055,bry:quillen}.

For any proper algebraic manifold $X$, the Hermitian holomorphic
Deligne cohomology groups of weight $p$ were originally
introduced by Brylinski in \cite{bry:quillen} as the
hypercohomology groups, denoted by $\dhhH{X}{p}$, of certain
complexes $\brydhh{X}{p}$, whose definition is reminiscent of the
corresponding one for the standard Deligne complexes. In
particular, if $\Pichat (X)$ denotes the group of isomorphism
classes of line bundles on $X$ equipped with a smooth hermitian
metric, it is shown that
\begin{equation*}
  \Pichat (X) \iso \dhhH[2]{X}{1}\,.
\end{equation*}
The groups $\dhhH{X}{p}$ have cup products behaving in the
standard way:
\begin{equation*}
  \dhhH[k]{X}{p} \otimes \dhhH[l]{X}{q}
  \xrightarrow{\cup} \dhhH[k+l]{X}{p+q}\,.
\end{equation*}
Thus, given the two metrized line bundles as above we can
multiply their classes in $\Pichat (X)$ using this product to
obtain a class $[\Bar{\sheaf{L}}]\cup [\Bar{\sheaf{M}}] \in
\dhhH[4]{X}{2}$.

If $X$ is a curve, the result follows since the latter group is
isomorphic to $\H^2(X, \RR)$. (See ref.\ \cite{math.CV/0211055}
and sect.~\ref{sec:prod-herm-line}, below.  We have neglected the
Tate twists by $\tate$, here.) The number so obtained is the
length of a generator of $\dual{\sheaf{L}}{\sheaf{M}}$ associated
to a specific choice of $l$ and $m$, see
sect.~\ref{sec:determ-cohom} below.

The case $\sheaf{L} =\sheaf{M} =T_X$ turns out to be quite
interesting on its own and relevant to the hyperbolic geometry of
$X$ considered as a Riemann surface.  A hermitian fiber metric is
now just a conformal metric on $X$. If the genus is $\geq 2$, we
have shown in \cite[Thm.  5.1]{math.CV/0211055} that the cup
square $[T_X]\cup [T_X]$, supplemented by the area form,
determines a functional whose extremum is precisely the
hyperbolic metric of constant curvature equal to $-1$. In fact,
we have shown it coincides with the Liouville functional studied
in refs.\ \cite{zogtak1987-1,zogtak1987-2,math.CV/0204318}, in
relation to Weil-Petersson geometry.

Returning to general case of the pairing
$\dual{\sheaf{L}}{\sheaf{M}}$ where $\sheaf{L}$ and $\sheaf{M}$
are not necessary equal, it has been observed---also by Deligne,
\cite{MR89b:32038}---that the metric on the pairing could be
defined also when the fiber metrics on both $\sheaf{L}$ and
$\sheaf{M}$ are allowed to be singular, provided the
corresponding loci are disjoint. In a different direction, if
$\mathscr{X}/_\ZZ$ is an arithmetic surface such that $X$ is the
corresponding fiber at infinity, U. K\"uhn (\cite{MR2002d:14035})
has formulated a pairing at arithmetic infinity (generalizing
Deligne's one) for line bundles where the metric has logarithmic
singularities at a finite set of points.  (This is shown in the
same work to be compatible with an earlier version by Bost, cf.\
\cite{MR2000c:14033}, based on Green's currents.)

This leads to the interesting question of whether pairings with
singular metrics can be expressed in terms of a natural cup
product, and if the results outlined above extend to singular
metrics. We show this is indeed the case.

In this work we first present an extension of the framework of
Hermitian holormorphic Deligne cohomology to pairs $(X,D)$, where
$X$ is a complete complex algebraic manifold, and $D$ is a normal
crossing divisor. Then, by working on a regular projective curve
over \CC, we use the cup product in this cohomology to obtain the
pairing for two hermitian line bundles with hermitian metrics
singular along $D$.

In slightly more detail, we define a new model $\ndhh{X}{p}$ for
the Hermitian Deligne complex (namely, we replace the complex
$\brydhh{X}{p}$ above by a quasi-isomorphic one) and obtain the
corresponding cup product structure in Theorem~\ref{thm:1}. In
fact, this is first defined for the case where $D$ is empty, a
case which is already of independent interest. In this case
$\ndhh{X}{p}$ is simply the cone of a morphism between two
complexes,
\begin{equation*}
  \ndhh{X}{p} =
  \cone \bigl(
  \delz{p} \lto \sigma^{<2p}\delalg{\sha{X}}{p}
  \bigr)[-1] \,,
\end{equation*}
where $\delz{p}$ is the \ZZ-valued Deligne complex, and
$\sigma^{<2p}\delalg{\sha{X}}{p}$ is a suitable truncation of a
Deligne algebra based on the Dolbeault complex of $X$. When
$D\neq \varnothing$, the Hermitian Deligne complex and the
corresponding (relative) cohomologies are immediately defined for
the pair $(X, U=X\setminus D)$ in terms of a sheaf
$\ndhh{X,U}{p}$ on a pair of spaces in the sense of ref.\
\cite{MR86h:11103}.

A further refinement is obtained by imposing growth conditions by
considering appropriate subsheaves of $\jmath_*
\sigma^{<2p}\delalg{\sha{U}}{p}$. We are ultimately interested in
good metrics in the sense of Mumford (\cite{MR81a:32026}), whose
characteristic forms have singularities akin to those of the
standard Poincar\'e metric on a punctured disk. While we are able
to directly use good forms in complex dimension $1$, it is not
possible to do so to set the homological algebra machinery in
general. In order to proceed, we exploit the complex of
\emph{pre-log-log} forms recently introduced in ref.\
\cite{math.AG/0404122}. It is possible to consider slightly
different kinds of growth conditions: in particular, the complex
of \emph{log-log} forms, also mentioned in loc.\ cit., seems to
have several advantages, including a Poincar\'e
lemma.\footnote{\label{fn:1} A published definition is not
  available at the time of this writing. A paper by the same
  authors of ref.\ \cite{math.AG/0404122} is in preparation. We
  warmly thank Prof.  Burgos for kindly informing us.}  The issue
is briefly touched upon in sect.~\ref{sec:product-curve}, but we
have not dwelt on it, since we were able to work with good forms
to a greater extent than the considerably more difficult ref.\
\cite{math.AG/0404122}.

Having set up the general framework of Hermitian holomorphic
Deligne cohomology with growth condition in general in
sect.~\ref{sec:sing-herm-struct}, we state our main applications
and results in the case of a complete regular complex curve $X$,
where $D$ reduces to a finite set of points.  If now $\sheaf{L}$
and $\sheaf{M}$ are line bundles on $X$ carrying a good metric
along $D$, we first show in Theorem.~\ref{thm:2} that the cup
product in Hermitian Deligne cohomology extends the Deligne
pairing to this case, thus generalizing the corresponding result
obtained in ref.\ \cite{math.CV/0211055} for the empty divisor.
As a consequence, we can obtain a generalization of the earlier
result in ref.\cite{math.CV/0211055} concerning hyperbolic
metrics.  Namely, good metrics $\d s^2$ on $T_U$ correspond to
the extension $\sheaf{L}=T_X(-D)$.  We obtain in
Theorem~\ref{thm:3} that the cup square $[\sheaf{L}] \cup
[\sheaf{L}]$ supplemented by the area integral (which is well
defined, since good forms are locally integrable) is still the
generating functional for the uniformizing metric. We do not need
the correction terms required in earlier definitions (cf. e.g.
\cite{zogtak1987-1}).

\subsection{Organization}
\label{sec:organization}

Let us briefly illustrate the organization of this work. The
definition of the Hermitian holomorphic Deligne complex
$\ndhh{X}{p}$ and the corresponding cohomology groups, as well as
the calculation of the multiplicative structure,  are carried
out in section~\ref{sec:models-herm-holom}, where we state
Theorem~\ref{thm:1}. We only sketch the proof, since many of the
arguments are either already available in the literature or
straightforward, but long to reproduce in detail.

Definitions are initially given in the proper case, without
reference to the divisor $D$. This is remedied in
section~\ref{sec:sing-herm-struct}, where definitions are
extended to pairs $(X,D)$. We also recall the relevant notions
necessary to consider specific growth conditions: Poincar\'e
forms, good forms, good metrics, pre-log-log forms.

Section~\ref{sec:prod-herm-line} collects some intermediate
results to be used later. We specialize the framework developed
in sections~\ref{sec:models-herm-holom}
and~\ref{sec:sing-herm-struct} to calculate the cup product
$[\Bar{\sheaf{L}}] \cup [\Bar{\sheaf{M}}]$ of two good hermitian
line bundles $\Bar{\sheaf{L}} = (\sheaf{L}, \norm{\phantom{s}})$
and $\Bar{\sheaf{M}} = (\sheaf{M}, \norm{\phantom{s}})$, both on
a general $(X,D)$ and then further specializing to the case $\dim
X=1$. We also collect some remarks about the integration map in
dimension $1$ which is also needed in the subsequent sections.

Sections~\ref{sec:deligne-pairing} and~\ref{sec:extr-hyperb-metr}
are devoted to our main results. In the former we extend the
Deligne pairing to good hermitian line bundles, and in the latter
we show that, after correcting with an area term, the pairing is
the generating functional for the constant negative curvature
hyperbolic metric. By then, most of the necessary preliminaries
will have already been obtained, therefore the corresponding
proofs have been kept to a minimum.

\subsection{Conventions and notations}
\label{sec:conv-notat}

For a subring $A$ of $\RR$ and an integer $p$, $A(p) =
(\tate)^p\,A$ is the Tate twist of $A$. We identify $\CC/\ZZ(p)
\iso \CC^\unit$ via the exponential map $z \mapsto \exp
(z/(\tate)^{p-1})$, and $\CC\iso \RR(p) \oplus \RR (p-1)$, so
$\CC / \RR(p) \iso \RR (p-1)$.

The projection $\pi_p\colon \CC \to \RR(p)$ is given by
$\pi\sb{p} (z) = \onehalf ( z + (-1)\sp p \Bar z)$, for $z\in
\CC$, and similarly for any other complex quantity, e.g.  complex
valued differential forms.

If $X$ is a complex manifold, $\shomega{X}$ denotes the de~Rham
complex of holomorphic forms, where we set $\sho{X} \coin
\shomega[0]{X}$ as usual.  $\she{X}$ denotes the de~Rham complex
of sheaves of $\RR$-valued smooth forms on the underlying smooth
manifold. Furthermore, $\sha{X} = \she{X}\otimes_\RR \CC$, and is
$\she{X} (p)$ the twist $\she{X} \otimes_\RR \RR(p)$.  Also,
$\sha[p,q]{X}$ will denote the sheaf of smooth $(p,q)$-forms, and
$\sha[n]{X} = \bigoplus_{p+q=n} \sha[p,q]{X}$, where the
differential decomposes in the standard fashion, $\d=\del +
\delb$, according to types. We also introduce the imaginary
operator $\dc = \del -\delb$ (with a slight departure from
convention, we omit the customary factor $1/(4\pi \sqrt{-1})$).
We have the rules
\begin{displaymath}
  \d\pi_p(\omega) = \pi_p(\d\omega)\,,\quad
  \dc\pi_p(\omega) = \pi_{p+1}(\dc\omega)
\end{displaymath}
for any complex form $\omega$.  Note that we have $2\del\delb =
\dc\d$.

We denote the complexes of global sections arising from the above
mentioned sheaf complexes by corresponding straight letters, as
$\complex{A}(X)$, $\complex[p,q]{A}(X)$, $\complex{E}(X)$,
$\complex{E}(X)(p) \eqdef \complex{E}(X)\otimes_\RR \RR(p)$, and
so on.

An open cover of $X$ will be denoted by $\cover{U}_X$. If
$\{U_i\}_{i\in I}$ is the corresponding collection of open sets,
we write $U_{ij} = U_i\cap U_j$, $U_{ijk} = U_i\cap U_j\cap
U_k$, and so on. More generally we can also have $\cover{U}_X =
\{ U_i \to X\}_{i\in I}$, where the maps are regular coverings in
an appropriate category. In this case intersections are replaced
by $(n+1)$-fold fibered products
\begin{math}
  U_{i_0 i_1\dotsb i_n}
  = U_{i_0} \times_X\dotsb \times_X U_{i_n}\,.
\end{math}

If $\complex{\sheaf{F}}$ is a complex of abelian sheaves on $X$, with
differential $d_\sheaf{F}$,
its \cech\ resolution with respect to a covering $\cover{U}_X\to
X$ is the double complex
\begin{displaymath}
  \CCC[{p,q}] (\sheaf{F}) \eqdef
  \vC[q]{\cover{U}_X}{\sheaf{F}^p}\,,
\end{displaymath}
where the $q$-cochains with values in $\sheaf{F}^p$ are given by
\begin{math}
  \prod \sheaf{F}^p (U_{i_0\dotsb i_n})\,.
\end{math}
The \cech\ coboundary operator is denoted $\deltacheck$. The
convention we use is to put the index along the \cech\ resolution
in the \emph{second} place, so the total differential is given by
$D=d_\sheaf{F} + (-1)^p \deltacheck$ on the component
$\vC[q]{\cover{U}_X}{\sheaf{F}^p}$ of the total simple complex.
Because of the Koszul sign rule we get the following convention
for \cech resolutions of complexes of sheaves.  If
$\sheaf{G}^\bullet$ is a second complex of sheaves on $X$, then
one defines the cup product
\begin{displaymath}
  \cup : \CCC[p,q](\sheaf{F}) \otimes \CCC[r,s](\sheaf{G})
  \lto
  \vC[q+s]{\cover{U}_X}{\sheaf{F}^p\otimes \sheaf{G}^r} \subset
  \CCC[p+r,q+s](\sheaf{F}\otimes\sheaf{G})
\end{displaymath}
of two elements $\{f_{i_0,\dotsc,i_q}\}\in \CCC[p,q](\sheaf{F})$
and $\{g_{j_0,\dotsc,j_s}\} \in \CCC[r,s](\sheaf{G})$ by
\begin{displaymath}
  (-1)^{qr}\,f_{i_0,\dots,i_q}\otimes
  g_{i_q,i_{q+1},\dotsc,i_{q+s}} \,.
\end{displaymath}
For a given complex of abelian objects, say $\complex{A}$, we
have the sharp truncations $\sigma^p\coin \sigma^{\geq p}$ and
$\sigma^{<p}$, namely $\sigma^p\complex[n]{A}=0$ for $n<p$, and
$\sigma^{<p}\complex[n]{A}=0$ if $n\geq p$.  On occasion we will
use the short-hand notations $\complex[\bullet\geq p]{A}$ and
$\complex[\bullet< p]{A}$, respectively.

Finally, we denote a quasi-isomorphism with $\lqi$, and a
homotopy equivalence by $\hoeq$.

\subsection*{Acknowledgments}
We thank the International School for Advanced Studies (SISSA) in
Trieste, Italy, for hospitality while part of this work was
written. We also thank Paolo Aluffi and Ugo Bruzzo for
discussions on the content of this paper. We would also like to
warmly thank J. Burgos Gil for an email exchange answering
several questions.

\section{Models of Hermitian-holomorphic Deligne cohomology}
\label{sec:models-herm-holom}

\subsection{Deligne cohomology}
\label{sec:deligne-cohomology}

In the following, $X$ will be a complete complex algebraic
manifold, that is the analytic variety associated to a smooth
proper scheme over \CC. (Since in the following we will work in
the analytic category there will be no need to distinguish
notationally between $X$ and $X_\mathit{an}$.)

First, recall that for a subring $A\subset \RR$ and an integer
$p$, the Deligne cohomology groups of weight $p$ of $X$ with
values in $A$ are the hypercohomology groups
\begin{equation}
  \label{eq:1}
  \delH{X}{A}{p} \eqdef \HHH^{\,\bullet} (X,\deligne{}{A}{p})\,,
\end{equation}
where $\deligne{}{A}{p}$ is the complex
\begin{equation}
  \label{eq:2}
  \deligne{}{A}{p} =
  \cone \bigl( A(p)_X\oplus F^p\shomega{X} \lto \shomega{X} \bigr)[-1]
\end{equation}
(the map in the cone is the difference of the two inclusions) and
$F^p\shomega{X}\coin \sigma^p \shomega{X}$ is the Hodge (``b\^ete'')
filtration. In the derived category the complex in \eqref{eq:2} is
isomorphic to:
\begin{equation}
  \label{eq:3}
  0\lto A(p)_X
  \overset{\imath}{\lto} \sho{X}
  \overset{\d}{\lto} \shomega[1]{X}
  \overset{\d}{\lto} \dotsm
  \overset{\d}{\lto} \shomega[{p-1}]{X}\,.
\end{equation}
Deligne cohomology has a graded commutative cup product
\begin{equation}
  \label{eq:4}
  \delH[k]{X}{A}{p}\otimes \delH[l]{X}{A}{q}
  \overset{\cup}{\lto} 
  \delH[k+l]{X}{A}{p+q}
\end{equation}
which is induced by a family $\alpha\rightsquigarrow\cup_\alpha$ of
cup products at the level of complexes
\begin{equation*}
  \deligne{}{A}{p}\otimes \deligne{}{A}{q}
  \overset{\cup_\alpha}{\lto} \deligne{}{A}{p+q}\,,
\end{equation*}
where $\alpha\in [0,1]$. This results from a general prescription
to obtain cup products on cones of certain diagrams of complexes
devised by \bei\cite{bei:hodge_coho}, and to be recalled below.
In particular the product $\cup_0$ induces the following cup
product on the (simpler) complex~\eqref{eq:3}. If $a\in
\deligne{}{A}{p}$ and $b\in \deligne{}{A}{q}$, then from
ref.~\cite{MR86h:11103,esn-vie:del} we have:
\begin{equation}
  \label{eq:5}
  a\cup b =
  \begin{cases}
    a\cdot b & \deg a = 0\,,\\
    a\wedge \d b & \deg a > 0\;
    \text{and}\; \deg b =q\,,\\
    0 &\text{otherwise.}
  \end{cases}
\end{equation}

\subsection{Hermitian holomorphic variants}
\label{sec:herm-holom-vari}

\subsubsection{}
\label{sec:pichat}

Let $\Pichat (X)$ be the group of isomorphism classes of line bundles
on $X$ with smooth hermitian fiber metric, that is of pairs
$(\sheaf{L},\rho)$, where $\sheaf{L}$ is an invertible
$\sho{X}$-module equipped with a map
\begin{math}
  \rho \colon \sheaf{L} \lto \she[+]{X}\,,
\end{math}
into (the sheaf of) positive real smooth functions.  A basic
observation by Deligne is that $\Pichat (X)$ can be identified with
the hypercohomology group:
\begin{equation}
  \label{eq:6}
  \HHH^2\bigl(X,
  \ZZ(1)_X \overset{\imath}{\lto} \sho{X}
  \xrightarrow{\pi_0} \she[0]{X}\bigr)\,,
\end{equation}
which is easy to see in \cech cohomology.  Indeed, for a cover
$\cover{U}_X$ of $X$, and a pair $(\sheaf{L},\rho)$, let $s_i$ be a
trivialization of $\sheaf{L}\rvert_{U_i}$, with transition functions
$g_{ij}\in \sho{X}^\unit (U_{ij})$ determined by $s_j = s_i g_{ij}$.
Let $\rho_i$ be the value of the quadratic form associated to $\rho$
on $s_i$, namely $\rho_i = \rho(s_i)$. Then we have $\rho_j = \rho_i
\abs{g_{ij}}^2$. Taking logarithms, we see that
\begin{equation}
  \label{eq:7}
  \bigl(\tate c_{ijk},
  \log g_{ij}, -\tfrac{1}{2} \log \rho_i \bigr) \,, 
\end{equation}
where
\begin{math}
  \tate c_{ijk} = \log g_{jk} -\log g_{ik} +\log g_{ij}\in
  \ZZ(1)\,,
\end{math}
is a cocycle representing the class of the pair $(\sheaf{L},\rho)$.

\subsubsection{}

Hermitian-holomorphic Deligne cohomology groups are a
generalization of the above fact to any weight and
degree. Namely, consider the complexes
\begin{equation}
  \label{eq:8}
  \brydhh{}{p} = \cone \bigl(
  \ZZ(p)_X \oplus (F^p\!\sha{X}\cap \sigma^{2p}\she{X}(p))
  \lto \she{X}(p)
  \bigr)[-1]\,,
\end{equation}
where $F^p\!\sha{X} = \bigoplus_{p'\geq p} \sha[p',\bullet-p']{X}$ is
the Hodge filtration.
\begin{definition}
  \label{def:1}
  The \emph{Hermitian-holomorphic Deligne} cohomology groups are the
  hypercohomology groups
  \begin{equation}
    \label{eq:9}
    \dhhH{X}{p} \eqdef \HHH^{\,\bullet} (X,\brydhh{}{p})
  \end{equation}
\end{definition}
The complexes~\eqref{eq:8} were introduced in ref.~\cite{bry:quillen}
by Brylinski.  It is easy to see that $\brydhh{}{1}$ is
quasi-isomorphic to the complex appearing in~\eqref{eq:6}, so that
$\Pichat (X)$ is identified with $\dhhH[2]{X}{1}$. More generally, the
elements of the higher degree groups $\dhhH[k]{X}{p}$ for $k=3$ (resp.
$k=4$) and $p=1,2$ correspond to classes of analytic (or algebraic)
gerbes (resp. $2$-gerbes) on $X$ with an appropriate notion of
hermitian and canonical connective structures (cf. ref.
\cite{math.CT/0310027}).

It is advantageous to have alternative ways of calculating the groups
$\dhhH{X}{p}$ by means of complexes quasi-isomorphic to
$\brydhh{}{p}$. In ref.~\cite{math.CV/0211055} we introduced the
complexes
\begin{equation}
  \label{eq:10}
  \dhh{}{p} = 
  \cone \bigl(
  \delz{p}\oplus (F^p\!\sha{X}\cap \sigma^{2p}\she{X}(p))
  \lto \deltilde{}{\RR}{p}
  \bigr)[-1]\,,
\end{equation}
where the complex
\begin{equation*}
  \deltilde{}{\RR}{p} =
  \cone \big(F^p\!\sha{X} \rightarrow \she{X}(p-1)\big)[-1]
\end{equation*}
also computes the real Deligne cohomology,
cf.~\cite{MR86h:11103,esn-vie:del}.  The complex $\dhh{}{1}$
provides a characterization of the \emph{canonical connection}
associated to the hermitian fiber metric structure. In fact the
quasi-isomorphism between $\dhh{}{1}$ and the complex
in~\eqref{eq:6} allows to conclude that the canonical connection
is unique. There are other advantages in using the complexes
$\dhh{}{p}$, discussed in detail in \cite{math.CV/0211055},
related to the fact that analytic or algebraic structures can be
described directly without resorting to the underlying smooth
objects.  This, however, happens at the cost of quite an amount
of algebraic intricacy, especially in the product structure.

Both complexes~\eqref{eq:8} and~\eqref{eq:10} admit cup products
according to \bei's general prescription, so that there results a
corresponding graded commutative product in cohomology:
\begin{equation}
  \label{eq:11}
  \dhhH[k]{X}{p}\otimes \dhhH[l]{X}{q}
  \overset{\cup}{\lto} 
  \dhhH[k+l]{X}{p+q}\,.
\end{equation}
We refer to refs. \cite{math.CV/0211055,math.CT/0310027} for the
details about the quasi-isomorphism $\brydhh{}{p}\lqi \dhh{}{p}$ and
more generally about the properties mentioned above.

\subsubsection{}

We want to introduce a third---and simpler---model for hermitian
holomorphic Deligne cohomology which is based on the following
observation.  Real Deligne cohomology can also be computed by the
following complex considered in detail in refs.
\cite{MR99d:14015,math.AG/0207036} (its use had been previously
suggested by Deligne):
\begin{equation}
  \label{eq:12}
  \delalg[n]{\sha{X}}{p} = 
  \begin{cases}
    \displaystyle
    \she[n-1]{X}(p-1)\bigcap
    \bigoplus_{\substack{p'+q'=n-1\\ p'< p, q'< p}}
    \sha[p',q']{X} & \text{if $n\leq 2p-1$} \\
    \displaystyle
    \she[n]{X}(p) \bigcap
    \bigoplus_{\substack{p'+q'=n\\ p'\geq p, q'\geq p}}
    \sha[p',q']{X} & \text{if $n\geq 2p$}
  \end{cases}
\end{equation}
with differential
\begin{equation}
  \label{eq:13}
  d_\goth{D} = 
  \begin{cases}
    -\pi\circ \d & n < 2p-1\\
    -2\delb\del & n = 2p-1\\
    \d           & n = 2p\,,
  \end{cases}
\end{equation}
where in the first line of~\eqref{eq:13} $\pi$ is the
projection.\footnote{Some sign differences result from different
  conventions about cones.  Given $f:\complex{A}\to \complex{B}$,
  ref.  \cite{MR99d:14015} uses $\cone (-f)[-1]$, whereas we use
  $\cone (f)[-1]$.}  The quasi-isomorphism is the composition
\begin{equation*}
  \delr{p} \lqi
  \cone \bigl(\she{X}(p)\oplus F^p\!\sha{X} \lto \sha{X}
  \bigr)[-1]
  \lqi \delalg{\sha{X}}{p}\,.
\end{equation*}
The left map trivially follows from $\RR(p)\qi \she{X}(p)$,
\eqref{eq:2}, and the fact that the inclusion
$F^p\shomega{X}\hookrightarrow F^p \sha{X}$ is a filtered
quasi-isomorphism. The right one is in fact a homotopy equivalence. It
is computed in detail in ref. \cite{MR99d:14015}, where the homotopies
are also explicitly computed.

Thus we can replace $\deltilde{}{\RR}{p}$ with
$\delalg{\sha{X}}{p}$ in the definition~\eqref{eq:10}, and still
obtain the same cohomology groups.  From the obvious map
$\delz{p}\to \delr{p}$, and following the chain of
quasi-isomorphisms, one has the map
\begin{equation*}
  \rho_p \colon \delz{p} \lto \delalg{\sha{X}}{p}
\end{equation*}
given by:
\begin{equation}
  \label{eq:14}
  \rho_p^n = 
  \begin{cases}
    0 & \text{if $n=0$}\\
    (-1)^n\pi_{p-1} & \text{if $1\leq n\leq p$}\\
    0 & \text{if $n\geq p$,}\\
  \end{cases}
\end{equation}
where we use version~\eqref{eq:3} of the Deligne complex.

Note that the sharp truncation $\sigma^{<2p}\delalg{\sha{X}}{p}$
is just the first line in \eqref{eq:12}, which comprises forms of
degree up to $2p-2$, whereas the model~\eqref{eq:3} of the
Deligne complex is zero after degree $p-1$. Hence the map $\rho_p$
factors through $\sigma^{<2p}\delalg{\sha{X}}{p}$, so we actually
have a map
\begin{equation}
  \label{eq:15}
  \rho_p \colon \delz{p} \lto \sigma^{<2p} \delalg{\sha{X}}{p}
\end{equation}
which is obviously also given by~\eqref{eq:14}.

Consider the complex:
\begin{equation}
  \label{eq:16}
  \ndhh{}{p} \eqdef 
  \cone \bigl(
  \delz{p} \overset{\rho_p}{\lto} \sigma^{<2p}\delalg{\sha{X}}{p}
  \bigr)[-1] \,.
\end{equation}
We have
\begin{lemma}
  The complex $\ndhh{}{p}$ is quasi-isomorphic to $\dhh{}{p}$.
\end{lemma}
\begin{proof}
  One need only observe that
  \begin{equation*}
    \sigma^{2p}\delalg{\sha{X}}{p} \equiv
    \she[n]{X}(p) \bigcap
    \bigoplus_{\substack{p'+q'=n\\ p'\geq p, q'\geq p}}
    \sha[p',q']{X} = F^p\!\sha{X}\cap \sigma^{2p}\she{X}(p)\,,
  \end{equation*}
  so that, writing:
  \begin{equation*}
    \dhh{}{p} = 
    \cone \bigl( \delz{p} \lto
    \cone \bigl(
    F^p\!\sha{X}\cap \sigma^{2p}\she{X}(p)
    \lto \deltilde{}{\RR}{p} \bigr)
    \bigr)[-1]\,,  
  \end{equation*}
  and replacing $\deltilde{}{\RR}{p}$ with $\delalg{\sha{X}}{p}$,
  the inner cone above is
  \begin{equation*}
    \delalg{\sha{X}}{p} / (F^p\!\sha{X}\cap \sigma^{2p}\she{X}(p))
    \lqi \sigma^{<2p}\delalg{\sha{X}}{p}\,.
  \end{equation*}
\end{proof}
\noindent
It follows from the lemma that $\dhhH[k]{X}{p}=
\HHH^k(X,\ndhh{}{p})$. 
\begin{remark}
  \label{rem:1}
  Since in the complex $\ndhh{}{p}$ we have the truncation of
  $\delalg{X}{p}$ after degree $2p$, from the expression
  \eqref{eq:13} of the differential $d_\goth{D}$ we obtain a map
  \begin{equation}
    \label{eq:17}
    \ndhh{}{p} \xrightarrow{-2\delb\del} \sha[p,p]{X}[-2p]\cap
    \she[2p]{X}(p)_\mathit{cl}\,.
  \end{equation}
  This immediately gives a ``characteristic class''
  \begin{equation}
    \label{eq:18}
    \dhhH[2p]{X}{p} \lto A^{p,p}(X)_{\RR(p),\mathit{cl}}\,,
  \end{equation}
  into $\RR(p)$-valued \emph{closed} forms.
\end{remark}
\begin{remark}
  The truncation occurring in the complex~\eqref{eq:16} is
  similar (but not equal) to the ``truncated Deligne complex''
  $\delalg[\bullet\leq 2p]{\sha{X}}{p}$ considered by Goncharov
  for the Arakelov motivic complex in~\cite{math.AG/0207036}.
  
  Also observe that the complex~\eqref{eq:16} is different from
  the complex used to define the ``truncated relative homology''
  groups in refs.\ \cite{MR99d:14015,math.AG/0404122}, as in the
  latter the truncation occurs in the first complex in the cone.
\end{remark}
\begin{example}
  The weights $p=1,2$ will be of special interests later on.
  
  The complex $\ndhh{}{1}$ coincides precisely with the complex
  appearing in~\eqref{eq:6}:
  \begin{equation*}
    \ZZ(1)_X \overset{\imath}{\lto} \sho{X}
    \xrightarrow{\pi_0} \she[0]{X}\,.
  \end{equation*}
  The complex $\ndhh{}{2}$ is the cone (shifted by 1) of the map
  \begin{equation}
    \label{eq:22}
    \begin{CD}
      \ZZ(2)_X @>\imath>> \sho{X} @>\d>> \shomega[1]{X} \\
      & & @VV-\pi_1V @VV\pi_1V \\
      & & \she[0]{X}(1) @>-\d>> \she[1]{X}(1)
      @>-\pi\circ \d>> \she[2]{X}(1)\cap \sha[1,1]{X}
    \end{CD}
  \end{equation}
\end{example}

\subsection{Multiplicative structure}
\label{sec:mult-structure}

\subsubsection{}
The complex $\ndhh{}{p}$ has a relatively simple multiplicative
structure. In fact there exists a one-parameter family of product
structures which is inherited from that of $\dhh{}{p}$, thanks to
\bei's construction. We present here the (simpler) one that will
be used in the sequel. We use the model \eqref{eq:3} for the
Deligne complex $\delz{p}$ with the product \eqref{eq:5}. Our
notations are as follows: an element of degree $n$ in
$\ndhh[n]{}{p}$ is denoted by $(x_p,w_p)$, where $x_p\in
\delz[n]{p}$, and $w_p\in
\sigma^{<2p}\delalg[n-1]{\sha{X}}{p}$. Also, $w^{r,s}_p$ denotes
the $(r,s)$-component of $w_p$.
Details of the proof of the following theorem will appear
elsewhere.
\begin{theorem}
  \label{thm:1}
  There is a map of complexes
  \begin{equation*}
    \cuphat \colon
    \ndhh{}{p}\otimes \ndhh{}{q} \lto \ndhh{}{p+q}
  \end{equation*}
  which is homotopy graded commutative. Let $(x_p,w_p)\in
  \ndhh[n]{}{p}$ and $(x_q,w_q)\in \ndhh[m]{}{q}$. Then:
  \begin{align*}
    (x_p,0) \cuphat (x_q,0) &=
    \begin{cases}
      (x_p\cup x_q,0) &  \text{$n=0$ or $m=0$}\\
      (x_p\cup x_q, (-1)^m 2\pi_p(x_p) \wedge \pi_{q-1}(x_q)) &
      \text{if $n$ and $m\neq 0$}
    \end{cases}\\
    (x_p,0) \cuphat (0,w_q) &=
    \begin{cases}
      0 & \text{$n=0$}\\
      (0,\pi_{p-1}(x_p)\wedge (\del w_q^{q-1,m-q-1}-\delb
      w_q^{m-q-1,q-1})) & \text{$n=1,\dots,p-1$}\\
      (0,\pi_{p-1}(x_p)\wedge (\del w_q^{q-1,m-q-1}
      -\delb w_q^{m-q-1,q-1})
      + (-1)^p\pi_p(\d x_p)\wedge w_q)
      & \text{$n=p$} \\
      0 & \text{$n\geq p+1$}
    \end{cases}\\
    (0,w_p)\cuphat (x_q,0) &= 0\\
    (0,w_p)\cuphat (0,w_q) &= 
    \begin{cases}
      0 & \text{$m\neq 2q$} \\
      (0, w_p\wedge 2\delb\del w_q) & \text{$m=2q$}
    \end{cases}
  \end{align*}
  The map $\cuphat$ induces the product structure~\eqref{eq:11}.
\end{theorem}
\begin{corollary}
  The map $\cuphat$ induces the product~\eqref{eq:11}.
\end{corollary}

\subsubsection{Sketch of the proof of Thm. \ref{thm:1}}
In order to get the product structure in Thm. \ref{thm:1} one
combines the following steps.

First we require a variant of the \bei product
(\cite{bei:hodge_coho}) for the cones of certain triangular
diagrams of complexes introduced in \cite{math.CV/0211055}:
namely, for $i=1,2,3$ consider diagrams
\begin{equation}
  \label{eq:19}
  \mathcal{D}_i \eqdef
  \complex{X}_i \overset{f_i}{\lto} \complex{Z}_i
  \overset{g_i}{\longleftarrow} \complex{Y}_i
\end{equation}
with product structures $\complex{X}_1\otimes \complex{X}_2
\xrightarrow{\cup} \complex{X}_3$, and similarly for
$\complex{Y}_i$, $\complex{Z}_i$, compatible with the $f_i,g_i$
up to homotopy, that is:
\begin{align*}
  f_3\circ \cup - \cup \circ (f_1\otimes f_2)
  &= d_Z h + h\, d_{X_1\otimes X_2}\,, \\
  g_3\circ \cup - \cup \circ (g_1\otimes g_2)
  &= d_Z k + k\, d_{Y_1\otimes Y_2}\,,
\end{align*}
with homotopies:
\begin{equation*}
  h \colon \complex{\bigl( X_1\otimes X_2 \bigr)}
  \lto \complex[\bullet -1]{Z}_3\:,\quad
  k \colon \complex{\bigl( Y_1\otimes Y_2 \bigr)}
  \lto \complex[\bullet -1]{Z}_3\,.
\end{equation*}
The following lemma is immediately verified
(\cite{math.CV/0211055}):
\begin{lemma}
  \label{lem:1}
  Let $\Gamma(\mathcal{D}_i) = \cone (\complex{X}_i\oplus
  \complex{Y}_i \xrightarrow{f_i-g_i} \complex{Z}_i)[-1]$,
  $i=1,2,3$. For
  \begin{math}
    (x_i,y_i,z_i) \in \complex{X}_i \oplus \complex{Y}_i \oplus
    \complex[\bullet -1]{Z}_i\,,\;i=1,2\,,
  \end{math}
  and $\alpha\in [0,1]$ the formula:
  \begin{equation}
    \label{eq:20}
    \begin{split}
      (x_1,y_1,z_1) \cup_\alpha (x_2,y_2,z_2) =
      \Big(&x_1\cup x_2, y_1\cup y_2, \\
      &(-1)^{\deg (x_1)}
      \big((1-\alpha )f_1(x_1) + \alpha g_1(y_1) \big) \cup z_2 \\
      &\quad +z_1\cup \big( \alpha f_2(x_2) +
      (1-\alpha)g_2(y_2)\big)\\
      &\qquad -h(x_1\otimes x_2) +k(y_1\otimes y_2) \Big)\,.
    \end{split}   
  \end{equation}
  defines a family of products
  \begin{equation*}
    \Gamma(\mathcal{D}_1)\otimes \Gamma(\mathcal{D}_2)
    \xrightarrow{\cup_\alpha}
    \Gamma(\mathcal{D}_3)\,.
  \end{equation*}
  These products are homotopic to one another, and graded
  commutative up to homotopy. The homotopy formula is the same as
  that found in ref.~\cite{bei:hodge_coho}.
\end{lemma}
As a further technical ingredient, the following proposition is used
to transfer product structures among homotopically equivalent
complexes. (It is probably well-known: a proof is included in ref.
\cite{MR99d:14015}.)
\begin{proposition}
  \label{prop:1}
  Let $\complex{X}$ and $\complex{Y}$ be two homotopically equivalent
  complexes, and let $\phi: \complex{X}\to \complex{Y}$ and $\psi :
  \complex{Y} \to \complex{X}$ be the homotopy equivalences. Assume
  $\complex{Y}$ has a product structure $\cup$. Then the position
  \begin{equation*}
    \psi (\phi(x_1)\cup \phi (x_2))
  \end{equation*}
  provides a product structure on $\complex{X}$ which is homotopy
  associative and graded commutative (up to homotopy) if so is
  the one on $\complex{Y}$.
\end{proposition}

Secondly, starting from the complex $\dhh{}{p}$ and using the
quasi-isomorphism
\begin{equation*}
  \delr{p} \lqi \delalg{\sha{X}}{p}
\end{equation*}
we apply the first part and lemma \ref{lem:1} to the diagram
\begin{equation}
  \label{eq:21}
  \mathcal{D}_{p} \eqdef
  \delz{p} \xrightarrow{\rho_p} \delalg{\sha{X}}{p}
  \xleftarrow{\mathit{incl}} F^p\!\sha{X}\cap \sigma^{2p}\she{X}(p)\,.
\end{equation}
The product structure on $\delz{p}$ is still~\eqref{eq:5}, and
that on $\delalg{\sha{X}}{p}$ can be found in ref.
\cite{MR99d:14015}.  In particular, it coincides with the
standard wedge product in degrees bigger or equal $2p$, which is
the product structure for the included subcomplex
$F^p\!\sha{X}\cap \sigma^{2p}\she{X}(p) =
\sigma^{2p}\delalg{\sha{X}}{p}$. The inclusion in the right part
of~\eqref{eq:21} is obviously strictly compatible with the
products, whereas the map $\rho_p$ is only compatible up to
homotopy since $\delalg{\sha{X}}{p}$ is only isomorphic to $\cone
\bigl(\she{X}(p)\oplus F^p\!\sha{X} \lto \sha{X} \bigr)[-1]$ up
to homotopy. (The latter complex carries the strict \bei
product.) For $x_p\in \delz[n]{p}$ and $x_q\in \delz[m]{q}$ the
homotopy $h$ is:
\begin{equation*}
  h(x_p\otimes x_q) = 
  \begin{cases}
    0 & \text{$n=0$ or $m=0$} \\
    (-1)^{m-1} 2\pi_p(x_p)\wedge \pi_{q-1}(x_q)& \text{otherwise.}
  \end{cases}
\end{equation*}

Finally, if in a diagram $\mathcal{D}$ of the form~\eqref{eq:19} the
map $g:\complex{Y}\to \complex{Z}$ is injective, there is a homotopy
equivalence
\begin{equation*}
  \Gamma (\mathcal{D})\hoeq
  \cone \bigl( X\overset{f}{\lto} \coker{g} \bigr)[-1]
\end{equation*}
we can use to transfer the product structure~\eqref{eq:20} from the
cone of the diagram $\mathcal{D}$ to the complex on the right by
making use of proposition~\ref{prop:1}. By applying this procedure to
$\Gamma (\mathcal{D}_p)$ in~\eqref{eq:21} we obtain the desired
product structure on the complex $\ndhh{}{p}$.

\section{Singular hermitian structures}
\label{sec:sing-herm-struct}

\subsection{Relative cohomology}
\label{sec:relat-cohom}

Recall that $X$ is a projective algebraic manifold. Let now $D$
be a normal crossing divisor in $X$ and let $U$ be the complement
of $D$ in $X$. Let $\jmath : U\to X$ be the inclusion. As usual,
an open coordinate subset $V$ with coordinates $(z_1,\dots,z_n)$
is \emph{adapted} to $D$ if the divisor $D$ is locally given by
the equation $z_1\dotsm z_k = 0$.

By viewing $X$ as a compactification of $U$, following ref.
\cite{MR86h:11103}, we consider the following sheaf complex on
the triple $(U,X,\jmath)$:
\begin{equation}
  \label{eq:24}
  \ndhh{X,U}{p} \eqdef
  \cone\bigl(
  \deligne{X}{\ZZ}{p} \xrightarrow{\rho_p}
  \jmath_*\delalg[\bullet < 2p]{\sha{U}}{p}
  \bigr)[-1]\,.
\end{equation}
Let us emphasize the $\ZZ$-valued Deligne complex is the one for
$X$, whereas the smooth part is the (direct image of the)
truncation of the $p$-th Deligne algebra of the open complement
$U$ of $D$.
\begin{definition}
  The hypercohomology groups
  \begin{equation}
    \label{eq:25}
    \dhhH{X,U}{p} \eqdef \HHH^\bullet (X,\ndhh{X,U}{p})
  \end{equation}
  are the Hermitian-holomorphic Deligne cohomology groups of the
  triple $(U,X,\jmath)$.
\end{definition}
\begin{remark}
  The elements of the group $\dhhH[2]{X,U}{1}$, correspond to
  pairs $(\sheaf{L},\rho)$, where $\sheaf{L}$ is a line bundle
  over $X$, and $\rho$ is a smooth hermitian metric on the
  restriction $\sheaf{L}\vert_U$.
\end{remark}

\subsubsection{}\label{sec:rel-coho-formal}
The characteristic morphisms~\eqref{eq:17} and~\eqref{eq:18}
generalize in this case to maps:
\begin{gather*}
  \ndhh{X,U}{p} \xrightarrow{-2\delb\del}
  \jmath_*\bigl( \sha[p,p]{U}[-2p]\cap
  \she[2p]{U}(p)\bigr)_\mathit{cl}\\
  \intertext{and}
  \dhhH[2p]{X,U}{p} \lto A^{p,p}(U)_{\RR(p),\mathit{cl}}\,.
\end{gather*}
{F}rom the cone in~\eqref{eq:24} and the fact that $\jmath$ is
affine we get the standard long exact cohomology sequence:
\begin{equation*}
  \dotsm \lto
  \HHH^{k-1}(U,\delalg[\bullet<2p]{\sha{U}}{p})
  \lto \dhhH[k]{X,U}{p} \lto \delH[k]{X}{\ZZ}{p} \lto
  \HHH^{k}(U,\delalg[\bullet<2p]{\sha{U}}{p}) \lto \dotsm
\end{equation*}
Moreover, the arguments in refs. \cite{MR86h:11103}, and
especially \cite[\textsection 4]{esn-vie:del}, show that there is
a product map
\begin{equation*}
  \cuphat: \ndhh{X,U}{p}\otimes \ndhh{X,U}{q} \lto
  \ndhh{X,U}{p+q} 
\end{equation*}
whose expression is still given by the formula in
Thm~\ref{thm:1}. As a consequence we obtain a cup product for the
relative cohomology groups:
\begin{equation*}
  \dhhH[k]{X,U}{p}\otimes \dhhH[l]{X,U}{q}
  \overset{\cup}{\lto} 
  \dhhH[k+l]{X,U}{p+q}\,.
\end{equation*}

\subsubsection{}
A variation on the previous theme is to consider \emph{two}
divisors with normal crossings $D$ and $D'$. Let us assume
that $D\cup D'$, and $D\cap D'$ are also divisors with
normal crossings. Let $U$ and $U'$ be their complements in
$X$, with inclusion maps $\jmath$ and $\jmath'$. Then
$\ndhh{X,U_1}{p}$ by restriction yields $\ndhh{X,U\cap
  U'}{p}$, and similarly for $\ndhh{X,U'}{q}$. Thus we get a
product
\begin{equation*}
  \dhhH[k]{X,U}{p}\otimes \dhhH[l]{X,U'}{q}
  \overset{\cup}{\lto} 
  \dhhH[k+l]{X,U\cap U'}{p+q}\,.
\end{equation*}

\subsubsection{}

In order to consider pairs $(\sheaf{L},\rho)$ where the hermitian
metric $\rho$ is required to have a prescribed behavior along
$D$, we will need to consider appropriate subcomplexes of
$\jmath_*\delalg[\bullet < 2p]{\sha{U}}{p}$ with specified growth
conditions.
\begin{definition}
  \label{def:2}
  For $p\geq 0$ consider a subcomplex
  \begin{equation*}
    \subdelalg{X}{p}{?} \subset \jmath_*\delalg{\sha{U}}{p}
  \end{equation*}
  where the ``$?$'' in the subscript denotes a growth condition
  to be described below, such that:
  \begin{enumerate}
  \item \label{item:1} The morphism $\rho_p : \deligne{X}{\ZZ}{p}
    \to
    \jmath_*\delalg[\bullet < 2p]{\sha{U}}{p}$ factors through\\
    $\subdelalg[\bullet < 2p]{X}{p}{?}$;
  \item \label{item:2} The cup-product restricts to
    $\subdelalg{X}{p}{?}$, namely
    \begin{equation*}
      \begin{CD}
        \subdelalg{X}{p}{?}\bigotimes \subdelalg{X}{q}{?} &
        @>>> & \subdelalg{X}{p+q}{?}\\
        @VVV & & @VVV\\
        \jmath_*\delalg{\sha{U}}{p}\bigotimes
        \jmath_*\delalg{\sha{U}}{q} & @>>> &
        \jmath_*\delalg{\sha{U}}{p+q}\,;
      \end{CD}
    \end{equation*}
  \item \label{item:3}
    \begin{math}\displaystyle
      \subdelalg{X}{p}{?} = \jmath_*\delalg{\sha{U}}{p} \bigcap
      \lbrace \text{forms with specified growth along $D$}\rbrace
    \end{math}
  \end{enumerate}
  The \emph{Hermitian-holomorphic Deligne complex with growth} is
  the complex
  \begin{equation}
    \label{eq:30}
    \ndhh{X,?}{p} \eqdef
    \cone \bigl( \deligne{X}{\ZZ}{p} \lto
    \subdelalg[\bullet < 2p]{X}{p}{?} \bigr)[-1]
  \end{equation}
\end{definition}
One accordingly defines cohomology groups
\begin{equation*}
  \dhhH{X}{p}_? \eqdef \HHH^\bullet (X,\ndhh{X,?}{p})\,.
\end{equation*}
By restriction from the sheaf~\eqref{eq:24} defined on $(X,U)$
the groups $\dhhH{X}{p}_?$ will satisfy the same formal
properties as those spelled in sect.~\ref{sec:rel-coho-formal}.

\subsection{Good metrics}
\label{sec:good-metrics}

We need to recall Mumford's notion of \emph{good hermitian
  metric} \cite[\S 1]{MR81a:32026} in the case of a line bundle.

Let $V$ an adapted neighborhood of a point $x\in D$. Consider a
polycylinder neighborhood $P$ of $D$ in $V$, so that
\begin{equation*}
  P\cap U = (\Dot\Delta_a)^k \times (\Delta_a)^{n-k}\,,
\end{equation*}
where $\Delta$ (resp. $\Dot\Delta$) is a standard disk (resp.
punctured disk) of radius $r=\abs{z}<a$, centered at $z=0$.

On the punctured disk $\Dot\Delta_a$ we have the Poincar\'e
metric
\begin{equation*}
  d s^2 = \frac{\abs{\d z}^2}{\bigl(\abs{z}\log
  \abs{z}\bigr)^2}\,, 
\end{equation*}
and the disk $\Delta_a$ is equipped with the standard Euclidean
metric $\abs{\d z}^2$. Denote by $h_P$ the resulting product
metric on the punctured polycylinder $(\Dot\Delta_a)^k \times
(\Delta_a)^{n-k}$.
\begin{definition}
  \label{def:3}
  A smooth $q$-form $\omega$ on $U$ is said to have
  \emph{Poincar\'e growth along $D$} if there is a cover
  $\cover{U}_X$ of $X$ with adapted neighborhoods and
  polycylinders as above such that the following estimate holds:
  \begin{equation}
    \label{eq:27}
    \abs{\omega (\xi_1,\dots,\xi_q)}^2 \leq
    C\, h_P(\xi_1,\xi_1) \dotsm h_P(\xi_q,\xi_q)\,.
  \end{equation}
  ($C$ is a positive constant, and $\xi_1,\dots,\xi_q$ are
  tangent vectors at some point of $P\cap U$.)  A form $\omega$
  is a \emph{good form (along $D$)} if both $\omega$ and
  $\d\omega$ have Poincar\'e growth.
\end{definition}
Clearly, good forms form a differential graded algebra.  It is
proved in ref.  \cite{MR81a:32026} that they are locally
integrable and the associated currents have no residue.

Let now $\sheaf{L}$ be a line bundle on $X$ equipped with a
smooth hermitian metric $\rho$ on the restriction
$\sheaf{L}\vert_U$.
\begin{definition}
  \label{def:6}
  The metric $\rho$ is \emph{good along $D$} if for any adapted
  neighborhood $V$ with polycylinder $P$ as above, there is a
  non-vanishing section $s\in \sheaf{L}\vert_{P\cap U}$ such that:
  \begin{enumerate}
  \item $\norm{s}$, $\norm{s}^{-1} \leq C \,
    \bigl(\sum_{i=1}^k\log \abs{z_i}\bigr)^N$, for some $C>0$,
    $N\in \NN$.
  \item The forms $\del\log \norm{s}$ and $\delb\del\log
    \norm{s}$ are good forms along $D$.
  \end{enumerate}
\end{definition}
We have set $\norm{s}^2=\rho (s,s)$.  We will denote such a pair
$(\sheaf{L},\rho)$ by $\Bar{\sheaf{L}}$ and call it a \emph{good
  hermitian line bundle.}

Given the pair $(\sheaf{L}\vert_U,\rho)$, there is a unique
extension $\sheaf{L}$ to $X$ such that $\rho$ is good. Moreover,
if $s$ is a section of $\sheaf{L}$ on $X$, the current associated
to the good form $\delb\del\log \norm{s}^2$ represents the first
Chern class of $\sheaf{L}$. (For both facts, cf. ref. \cite[\S
1]{MR81a:32026}.) Owing to this fact, if $\Bar{\sheaf{L}}$ is a
good hermitian line bundle, we sometimes abuse the language and
denote the form $\delb\del\log \norm{s}^2$ by $c_1(\sheaf{L})$.
\begin{example}
  \label{ex:hyp-curves}
  Let $X$ be a smooth projective algebraic curve over \CC\ (a
  Riemann surface). By our standing hypotheses $D$ will
  correspond to a finite set of points on $X$. An adapted
  neighborhood $V$ of a point $p\in D$ is simply a disk
  $\Delta_a$ of radius $a$ with a coordinate function $z$ such
  that $z(p)=0$. Also, $V\cap U=\Dot\Delta_a$, the corresponding
  punctured disk.
  
  Consider the hyperbolic metric $\d s^2_\hyp$ on $T_U$ of
  constant negative curvature equal to $-1$.  Locally near $p\in
  D$ we have:
  \begin{equation*}
    \d s^2_\hyp\Bigr\vert_{V\cap U} = \rho_\hyp \abs{\d z}^2\,, 
  \end{equation*}
  where $\rho_\hyp$ is a smooth positive function on $V\cap U$,
  such that there exists a continuous function $\alpha:V\to \RR$,
  smooth on $V\cap U=\Dot\Delta_a$, and the following asymptotic
  behavior at $p$ (cf.\ \cite{zogtak1987-1}):
  \begin{equation}
    \label{eq:28}
    \log\rho_\hyp (z) = -\log
    \bigl( \abs{z}^2 \log^2 \abs{z} \bigr)+\alpha (z)\,,
    \quad
    \abs{\del \alpha}^2 \leq
    \frac{\abs{\d z}^2}{\abs{z}^2\log^4\abs{z}}\,.
  \end{equation}
  This is a good metric.  It is well-known that the good
  hermitian line bundle $\sheaf{L}$ on $X$ corresponding to the
  pair $(T_U,\d s^2_\hyp)$ is $\sheaf{L}=T_X(-D)$. Indeed, in the
  adapted neighborhood $V$, $s=z\frac{\del}{\del z}$ is a
  non-vanishing holomorphic section of $\sheaf{L}\vert_V$, and its
  length square with respect to the hyperbolic metric is
  $\norm{s}^2= \exp {\alpha} / \log^2 \abs{z}$, which implies
  goodness (cf.\ \cite{MR2002k:32016}).
\end{example}

\subsubsection{}\label{sec:not-good}
It would seem natural to use good forms to realize the complex of
Definition~\ref{def:2}.  Unfortunately, this is not possible, as
we now briefly explain.

Denote by $\sha{X}\lAngle D\rAngle_\good$ the
differential graded subalgebra of $\jmath_*\sha{U}$ comprised of
good forms.

If $\omega$ is a good form, it does not necessarily follow that
$\d\dc \omega=2\delb\del\omega$ is good too, namely it does not
necessarily have Poincar\'e growth along $D$. Hence the graded
module $\jmath_*\delalg{\sha{U}}{p}\bigcap \sha{X}\lAngle
D\rAngle_\good$ is not a complex, in general.  Note that
the form of the differential at degree $2p$ is the sole
troublesome point. Indeed, we have
\begin{lemma}
  \label{lem:2}
  The graded modules
  \begin{equation*}
    \jmath_*\delalg[\bullet<2p]{\sha{U}}{p} \bigcap
    \sha{X}\lAngle D\rAngle_\good\,,
    \quad
    \jmath_*\delalg[\bullet\geq 2p]{\sha{U}}{p} \bigcap
    \sha{X}\lAngle D\rAngle_\good
  \end{equation*}
  form well defined complexes.
\end{lemma}
\begin{proof}
  For the second, we have from eq.~\eqref{eq:13} that
  $d_\goth{D}=\d$, hence the notion of ``good'' form is
  compatible with the definition of the complex
  $\delalg{\sha{U}}{p}$ in degrees $\geq 2p$.
  
  For the first, let $\eta\in \jmath_*\sigma^{< 2p}
  \delalg[n]{\sha{U}}{p}$ be a good form. From~\eqref{eq:13}
  and~\eqref{eq:12} we see that to compute $d_\goth{D}\eta
  =-\pi(\d\eta)$ we must drop the component $\del\eta^{p-1,n-p}$ of
  $\d\eta$.  It is immediate to see that if $\d\eta$ has
  Poincar\'e growth, so does $\pi\circ\d\eta$.
\end{proof}
Moreover, the morphism $\rho_p\colon \deligne{X}{\ZZ}{p} \to
\jmath_*\delalg[\bullet < 2p]{\sha{U}}{p}$ clearly factors
through the subcomplex
\begin{equation*}
  \jmath_*\delalg[\bullet<2p]{\sha{U}}{p}
  \bigcap \sha{X}\lAngle D\rAngle_\good\,,
\end{equation*}
so that it is indeed possible to define the complex:
\begin{equation*}
  \ndhh{X,\good}{p} \eqdef
  \cone \bigl(
  \deligne{X}{\ZZ}{p} \lto
  \jmath_*\delalg[\bullet<2p]{\sha{U}}{p} \bigcap
  \sha{X}\lAngle D\rAngle_\good 
  \bigr)[-1]\,.
\end{equation*}
However, since in the expression of the product $\cuphat$ in
Theorem~\ref{thm:1} there appears the operator $2\delb\del$, we
see that condition~\ref{item:2} in Definition~\ref{def:2} is not
satisfied.

\subsection{Pre-log-log forms}
\label{sec:pre-log-log}

We consider a larger subcomplex (in fact a subalgebra) of
$\jmath_*\sha{U}$ of forms with a log-log-type growth which has
been recently introduced in ref. \cite{math.AG/0404122}. The
following two definitions are from loc.\ cit., \S 7, whose
notations we retain in part:
\begin{definition}
  \label{def:4}
  A smooth complex-valued function $f$ on $U$ has
  \emph{log-log-growth} along $D$ if
  \begin{equation*}
    f(z_1,\dots,z_n) \leq C\, \prod_{i=1}^k
    \log^N \abs{\log \abs{z_i}} 
  \end{equation*}
  on any adapted neighborhood $V$ with coordinates
  $z_1,\dots,z_n$ and some constants $C>0$ and $N\in \NN$.

  The \emph{sheaf of differential forms on $X$ with
  log-log-growth along $D$} is the subsheaf of $\jmath_*\sha{U}$
  generated by the functions with log-log-growth along $D$ and
  the differentials
  \begin{align*}
    & \frac{\d z_i}{z_i\,\log\abs{z_i}}\,,\;
    \frac{\d\Bar z_i}{\Bar z_i\,\log\abs{z_i}} & i&=1,\dots,k\,, \\
    & \d z_i\,,\; \d\Bar z_i & i&=k+1,\dots, n\,.
  \end{align*}
\end{definition}
\begin{definition}
  \label{def:5}
  A form $\omega$ with log-log-growth such that $\del\omega$,
  $\delb\omega$, and $\del\delb\omega$ also have log-log-growth
  along $D$ will be called a \emph{pre-log-log form.} The sheaf
  of pre-log-log forms is the subalgebra of $\jmath_*\sha{U}$
  generated by the pre-log-log forms. This complex is denoted by
  $\sha{X}\lAngle D\rAngle_\pre$. Similarly, the
  subcomplex of $\RR(p)$-valued forms will be denoted by
  $\she{X}(p)\lAngle D\rAngle_\pre$.
\end{definition}
We refer to op.\ cit.\ for more properties of $\sha{X}\lAngle
D\rAngle_\pre$. We only observe that while pre-log-log
forms and good forms are not the same, the former have some of
the salient features of the latter, notably, pre-log-log forms
are locally integrable and their associated currents do not have
residues (\cite[Proposition 7.6]{math.AG/0404122}).

Since $\sha{X}\lAngle D\rAngle_\pre$ has the same formal
properties of a Dolbeault complex, it makes sense to consider the
complex
\begin{math}
  \delalg{%
    \sha{X}\lAngle D\rAngle_\pre }{p}
\end{math}
which is a subalgebra of $\jmath_*\delalg{\sha{U}}{p}$ satisfying
the conditions in Definition~\ref{def:2}. Therefore, using the
same notations, we set:
\begin{equation}
  \label{eq:29}
  \subdelalg{X}{p}{\pre}\eqdef
  \delalg{\sha{X}\lAngle D\rAngle_\pre}{p}\,.
\end{equation}
\begin{definition}
  \label{def:7}
  The \emph{pre-log-log Hermitian-holomorphic Deligne complex} is
  the complex:
  \begin{equation}
    \label{eq:31}
    \ndhh{X,\pre}{p}
    = \cone \bigl(
    \deligne{X}{\ZZ}{p} \xrightarrow{\rho_p}
    \subdelalg[\bullet < 2p]{X}{p}{\pre} \bigr)[-1]\,.
  \end{equation}
  The hypercohomology groups:
  \begin{equation*}
    \dhhH{X}{p}_\pre \eqdef
    \HHH^\bullet (X,\ndhh{X,\pre}{p})
  \end{equation*}
  are the \emph{pre-log-log Hermitian-holomorphic Deligne
    cohomology groups.}
\end{definition}
In particular, good hermitian line bundles define elements of
pre-log-log Hermitian holomorphic Deligne groups:
\begin{proposition}
  \label{prop:2}
  Let $\Bar{\sheaf{L}}=(\sheaf{L},\rho)$ be a good hermitian line
  bundle.  Then $[\Bar{\sheaf{L}}]\in
  \dhhH[2]{X}{1}_\pre$.  Moreover, $c_1(\sheaf{L})\in
  \Gamma(X,\she{X}(1)\lAngle D\rAngle_\pre)$.
\end{proposition}
\begin{proof}
  Let $\cover{U}_X=\bigcup_i U_i$ be an open cover, and let
  $\sheaf{L}$ be trivialized by non-vanishing sections $s_i\in
  \sheaf{L}\vert_{U_i}$ as in sect.~\ref{sec:pichat} (with the
  same notations).
  
  It follows from definition~\ref{def:6} that since
  $\Bar{\sheaf{L}}$ is good, the $\log\rho_i$'s have log-log
  growth along $D\cap U_i$, and so do $\del\log\rho_i$ and
  $\del\delb\log\rho_i$. Hence $\log\rho_i$ is pre-log-log.
  
  Moreover, since
  \begin{equation*}
    c_1(\sheaf{L})\vert_{U_i} = \delb\del\log\rho_i
  \end{equation*}
  is of type $(1,1)$ and closed, good implies pre-log-log (cf.
  \cite{math.AG/0404122}).
\end{proof}

\section{Cup product of hermitian line bundles}
\label{sec:prod-herm-line}
Let us temporarily revert to the non-relative case.

We will need to have a closer look at the product
\begin{equation*}
  \Pichat (X)\coin \dhhH[2]{X}{1}
  \otimes \Pichat (X)\coin \dhhH[2]{X}{1}
  \lto \dhhH[4]{X}{2}\,.
\end{equation*}
induced by the cup product $\ndhh{}{1}\otimes \ndhh{}{1}
\xrightarrow{\cuphat} \ndhh{}{2}$, in particular we need an
explicit formula in \cech cohomology for the induced product at
the level of total \cech complexes.

Let $\Bar{\sheaf{L}} = (\sheaf{L},\rho)$ and $\Bar{\sheaf{L}'}
= (\sheaf{L}',\rho')$ be two hermitian line bundles on $X$, and
assume they are trivialized as in sect.  \ref{sec:pichat} with
respect to a cover $\cover{U}_X$.  Consider two \cech cocycles of
the form~\eqref{eq:7} representing the corresponding classes in
$\dhhH[2]{X}{1}$, where primed symbols refer to the second pair.

Using the product from Thm~\ref{thm:1}, and the conventions on
\cech resolutions outlined in sect.~\ref{sec:conv-notat}, the
cup-product $[\Bar{\sheaf{L}}]\cuphat [\Bar{\sheaf{L}'}]$ in
$\dhhH[4]{X}{2}$ is represented with a \cech resolution by the
cocycle:
\begin{equation}
  \label{eq:23}
  \renewcommand{\arraystretch}{1.3}
  \begin{array}{cc}
    (0,4) & (\tate)^2\,c_{ijk}\,c'_{klm} \\ \hline
    (1,3) & \tate\, c_{ijk} \log \, g'_{kl} \\ \hline
    (2,2) & (-\log \, g_{ij} \, \d\log\,g'_{jk}) \oplus
    \pi_1(\log\,g_{ij})\, \log\abs{g'_{jk}} \\ \hline
    (3,1) & -\log\abs{g_{ij}}\; \onehalf \dc\, \log\rho'_j
    +\dc\, \log\abs{g_{ij}} \onehalf \log\rho'_j
    \\ \hline 
    (4,0) & \onehalf \log\rho_i \; \delb\del\log\rho'_i
  \end{array}
\end{equation}
\begin{remark}
  \label{rem:2}
  From \eqref{eq:23} and remark~\ref{rem:1}, the characteristic form
  associated to $[\Bar{\sheaf{L}}] \cuphat [\Bar{\sheaf{L}'}]$ is the
  global form locally given by:
  \begin{equation}
    \label{eq:26}
    \delb\del\log\rho_i \wedge \delb\del\log\rho'_i\,.
  \end{equation}
  This form represents the product
  \begin{math}
    c_1(\sheaf{L}) \cup c_1 ( \sheaf{L}')
  \end{math}
  of the two Chern classes.  Moreover, according to refs.
  \cite{bry:quillen,math.CT/0310027}, this form is the Chern
  class of the hermitian $(2,2)$-curving on the $2$-Gerbe
  $\tame{\sheaf{L}}{\sheaf{L}'}$.
\end{remark}
Assume now $D$ is a normal crossing divisor in $X$ as before,
with $U=X\setminus D$, and let $\Bar{\sheaf{L}}$,
$\Bar{\sheaf{L}'}$ be good hermitian line bundles. The following
is evident:
\begin{proposition}
  \label{prop:3}
  The cup product
  \begin{equation*}
    \dhhH[2]{X}{1}_\pre \otimes
    \dhhH[2]{X}{1}_\pre \lto \dhhH[4]{X}{2}_\pre 
  \end{equation*}
  is computed again by eq.~\eqref{eq:23}, where elements in
  $\vC[4-k]{\cover{U}_X}{\ndhh[k]{X,\pre}{p}}$, for
  $k=3,4$, are pre-log-log forms along $D\cap U_J$, where $J$ is
  the multi-index $J=(i_0,i_1,\dots,i_{4-k})$.
\end{proposition}
In particular, combining this with proposition~\ref{prop:2}, we
immediately have that the cup product of good hermitian line
bundles is a pre-log-log hermitian holomorphic Deligne class,
whereas the characteristic form~\eqref{eq:26} is a pre-log-log
form on $X$ along $D$.
\begin{remark}
  It is immediately verified that the pre-log-log forms in the
  total cocycle~\eqref{eq:23} representing the class
  $[\Bar{\sheaf{L}}] \cuphat [\Bar{\sheaf{L}'}]$ as well as the
  characteristic form~\eqref{eq:26} are in fact good forms.
\end{remark}
In other words, given two good hermitian line bundles, their
cup-product can be represented in \cech cohomology by a total
cocycle with values in the complex $\ndhh{X,\good}{2}$
introduced in sect.~\ref{sec:not-good}.

\subsection{Product on a curve}
\label{sec:product-curve}

\subsubsection{}

Let us now consider in particular the case where $X$ is a smooth
proper curve over \CC\ as in example~\ref{ex:hyp-curves}. Since
$\dim_\RR X=2$, the complex $\ndhh{X,?}{2}$, where
$?=\good$ or $\pre$ (or any other growth
satisfying the requirements of Definition~\ref{def:2}), is the
cone shifted by 1 of:
\begin{equation*}
  \begin{CD}
    \ZZ(2)_X @>\imath>> \sho{X} @>\d>> \shomega[1]{X} \\
    & & @VV-\pi_1V @VV\pi_1V \\
    & & \she[0]{X}(1) \lAngle D\rAngle_? @>-\d>>
    \she[1]{X}(1) \lAngle D\rAngle_? @>- \d>>
    \she[2]{X}(1)\lAngle D\rAngle_?
  \end{CD}
\end{equation*}
where now we have the full complex $\she{X}(1)\lAngle
D\rAngle_\mathit{?}$. (So the choice $?=\good$ does make
sense in this case.)

Again by dimensional reasons, the part of the total
cocycle~\eqref{eq:23} corresponding to the element of
$\delH[4]{X}{\ZZ}{2}$ representing $\tame{\sheaf{L}}{\sheaf{L}'}$
becomes trivial, in fact we have:
\begin{lemma}
  \label{lem:3}
  The cocycle~\eqref{eq:23} reduces to a total \cech cocycle of
  degree $2$ in $\Tot\vC{\cover{U}_X}{\she{X}(1)\lAngle
    D\rAngle_?}$. That is, for a curve $X$ as above the cup
  product in Proposition~\ref{prop:3} reduces to:
  \begin{equation*}
    \dhhH[2]{X}{1}_?
    \otimes \dhhH[2]{X}{1}_?
    \xrightarrow{\cup}
    \HHH^2(X,\she{X}(1)\lAngle D\rAngle_?)
  \end{equation*}
\end{lemma}
\begin{proof}
  This follows immediately from a calculation identical to
  \cite[Thm. 5.1]{math.CV/0211055}, one need only observe that no
  problems with singularities at $D$ ever occur, since only
  sections of $\sho{X}$ and $\she[0]{X}(1)$ are involved.
  Namely, the cohomology class of the higher symbol
  $\tame{\sheaf{L}}{\sheaf{L}'}$ is zero and it can be written as
  the coboundary of a collection (of Bloch's type dilogarithms)
  $L_{ijk}\in \sho{X} (U_{ijk})$ such that:
  \begin{equation*}
    \d L_{ijk} = -\log g_{ij} \, \d\log g'_{jk}\,.
  \end{equation*}
  This allows to construct a cochain $\omega^2_{ijk}\in
  \she[0]{X}(1)(U_{ijk})\subset \she[0]{X}(1)\lAngle D\rAngle_?
  (U_{ijk})$ (in fact, a \cech cocycle) which, together with the
  $(4,0)$ and $(3,1)$ entries in~\eqref{eq:23}, respectively
  denoted $\omega^0_i$ and $\omega^1_{ij}$, forms a total degree
  $2$ cocycle $\Omega = \omega^0+\omega^1+\omega^2\in
  \Tot\vC{\cover{U}_X}{\she{X}(1)\lAngle D\rAngle_?}$. This
  cocycle can be injected via the standard cone exact sequence
  into $\Tot\vC{\cover{U}_X}{\ndhh{X,?}{2}}$: then its difference
  with~\eqref{eq:23} is shown to be a coboundary.
\end{proof}

\subsubsection{}

Depending on the chosen growth conditions, there remains the
question of whether in Lemma~\ref{lem:3} one has
\begin{equation*}
  \HHH^2(X,\she{X}(1)\lAngle D\rAngle_?) \iso \H^2(X,\RR(1))
\end{equation*}
(Of course, the latter is isomorphic to $\RR (1)$ if $X$ is
irreducible). This will be the case \emph{if} Poincar\'e lemma for
the complexes $\she{X}(p)\lAngle D\rAngle_?$ and $\sha{X}\lAngle
D\rAngle_?$ holds. In general, as observed in ref.\ \cite[remark
7.19]{math.AG/0404122}, $\H^2(X,\RR(1))$ must be at least a
direct factor of $\HHH^2(X,\she{X}(1)\lAngle D\rAngle_?)$, since
the composite map
\begin{equation*}
  \sha{X}\lto \sha{X}\lAngle D\rAngle_? \lto
  \complex{\sheaf{D}}_X \coin
  {}^\backprime\sha{X}[-2\dim X]
\end{equation*}
(where the complex on the right is the complex of complex-valued
currents) is a quasi-isomorphism. This remark can be generalized
to any $X$, not necessarily of complex dimension one, implying
that $\dhhH{X}{p}$ is a direct factor of $\dhhH{X}{p}_?$.

For pre-log-log forms, that is, if $?=\pre$, one indeed cannot
presently invoke Poincar\'e lemma. In ref.\
\cite{math.AG/0404122} it is observed this could be remedied by
resorting to a complex of forms with log-log-type growth
conditions imposed on all derivatives (cf.\
sect.\ref{sec:introduction}, footnote~\ref{fn:1}). We can work
around this issue with the integration map to be defined below.

\subsubsection{}
\label{sec:trace}

Returning to $X$ of complex dimension one, we can compose the
product map in Lemma~\ref{lem:3} with the map
\begin{math}
  \HHH^2(X,\she{X}(1) \lAngle D\rAngle_?) \lto \RR(1)
\end{math}
which can be defined as follows. First, observe that
\begin{equation*}
  \HHH^\bullet(X,\she{X}(1)\lAngle D\rAngle_?) \iso
  \H^\bullet (E^\bullet_X(1)\lAngle D\rAngle_?)\,,
\end{equation*}
where on the right we have the complex of global sections.  Let
$D=p_1+\dots +p_N$. Denote by $X_\varepsilon=X\setminus
B_\varepsilon (D)$, where $B_\varepsilon= \bigcup_{i=1}^N
\Delta_\varepsilon (p_i)$ is the union of disks of radius
$\varepsilon$ centered at each point $p_i$. If $\omega \in
E^2_X(1)\lAngle D\rAngle_?$ corresponds to a class in
$\HHH^2(X,\she{X}(1)\lAngle D\rAngle_?)$, then:
\begin{equation}
  \label{eq:42}
  \int_X \omega \eqdef \lim_{\varepsilon\to 0}
  \int_{X_\varepsilon} \omega\in \RR(1)\,.
\end{equation}
By local integrability of good or pre-log-log forms, this is well
defined.

Alternatively, let $\Sigma_\mathbf{T}$ be a representative of the
class $[X]$ obtained from a triangulation $\mathbf{T}$
subordinated to the cover $\cover{U}_X=\lbrace V_i\rbrace_{i\in I}$.
Thus there is a map of index sets
\begin{math}
  A \ni \alpha \mapsto i(\alpha)\in I
\end{math}
such that $\triangle^2_\alpha \subset V_{i(\alpha)}$,
$\triangle^1_{\alpha,\beta} \subset V_{i(\alpha) i(\beta)}$, etc.
Assume $\supp D\cap \mathbf{T}^{(1)} = \varnothing$, where
$\mathbf{T}^{(1)}$ is the $1$-skeleton of $\mathbf{T}$.

If $\Omega$ is the cocycle
representing a class $[\Omega]\in \HHH^{2}(X,\she{X}(1)\lAngle
D\rAngle_?)$, then set:
\begin{equation}
  \label{eq:43}
  \dual{\Omega}{\Sigma_\mathbf{T}}
  \eqdef
  \sum_\alpha \int_{\triangle^2_\alpha} \omega^0_{i(\alpha)}
  +\sum_{\langle\alpha\beta\rangle}
  \int_{\triangle^1_{\alpha\beta}} \omega^1_{i(\alpha)i(\beta)}
  +\sum_{\langle\alpha\beta\gamma\rangle}
  \int_{\triangle^0_{\alpha\beta\gamma}}
    \omega^2_{i(\alpha)i(\beta)i(\gamma)}\,.
\end{equation}
(Note that the sums and labels run over the abstract dual
triangulation.) More details can be found, e.g. in ref.
\cite{aldtak2000}.  A limiting process is implicit
in~\eqref{eq:43} as well. If $\triangle^2_\alpha$ contains $p\in
D$, then
\begin{equation}
  \label{eq:44}
  \int_{\triangle^2_\alpha} \omega^0_{i(\alpha)}
  \eqdef
  \lim_{\varepsilon\to 0} \int_{\triangle^2_\alpha\cap
  X_{\varepsilon}} \omega^0_{i(\alpha)}\,.
\end{equation}
In summary, we have:
\begin{lemma}
  \label{lem:5}
  Let $[\omega]=[\Omega]$. Then, with the above provisions:
  \begin{equation*}
    \int_X \omega = \dual{\Omega}{\Sigma_\mathbf{T}}\in \RR(1)\,.
  \end{equation*}
\end{lemma}
\begin{proof}
  The correspondence between $\Omega$ and $\omega$ can be made
  explicit using, e.g.  partitions of unity arguments. Then we
  need only observe that since $\mathbf{T}^{1}\cap \supp
  D=\varnothing$, singularities develop only in the integrations on
  the $2$-chains which hit points in $D$. These are handled
  by~\eqref{eq:44}.
\end{proof}

\section{Deligne pairing}
\label{sec:deligne-pairing}

\subsection{Reminder on Determinant of Cohomology}
\label{sec:determ-cohom}

$X$ is a smooth proper curve over \CC\ as in
example~\ref{ex:hyp-curves}. We need to collect a few formulas
and definitions from ref.\ \cite{MR89b:32038}.

Consider two line bundles $\sheaf{L}$ and $\sheaf{L'}$ on $X$ and
let $E$ and $E'$ be the corresponding divisors, assumed to have
disjoint support.  The Deligne pairing
$\dual{\sheaf{L}}{\sheaf{L'}}$, is a complex line generated by
the symbol $\dual{s}{s'}$, where $s$ (resp.\ $s'$) is a rational
section of $\sheaf{L}$ (resp.\ $\sheaf{L'}$), such that their
divisors are $(s)=E$ and $(s')=E'$. Replacing $s$ with $fs$,
where $f$ is a rational function, yields the relation
\begin{equation*}
  \dual{fs}{s'} = f(E') \; \dual{s}{s'}\,,
\end{equation*}
and similarly if we replace $s'$ by $gs'$. Consistency is ensured
by the Weil reciprocity property $f((g)) = g((f))$
\cite{gh:alg_geom}.

If $\sheaf{L}$ (resp. $\sheaf{L}'$) is equipped with a smooth
hermitian metric $\rho$ (resp. $\rho'$) over $X$, the complex
line $\dual{\sheaf{L}}{\sheaf{L}'}$ acquires a hermitian metric
given by (\cite{MR89b:32038}):
\begin{equation}
  \label{eq:32}
  \begin{split}
    \log\, \norm{\dual{s}{s'}}^2
    = \frac{1}{\tate} \int_X& \del\delb\log \norm{s}^2
    \,\log\norm{s'}^2\\ 
     &\qquad + \log\norm{s}^2[E'] +\log\norm{s'}^2[E]\,.
   \end{split}
\end{equation}
Since there is no danger of confusion, we have denoted all the
square norms corresponding to the various metrics by
$\norm{\phantom{s}}$. The operator $\del\delb$ is to be
interpreted in the sense of distributions.
 
Alternatively, let $f_1,f_2$ be two smooth $\RR_{\geq 0}$-valued
functions such that $f_1+f_2=1$ and $f_1$ (resp. $f_2$) vanishes
in a neighborhood of the support of $E'$ (resp. $E$). (It follows
that $f_1$ (resp. $f_2$) is equal to $1$ near the support of $E$
(resp. $E'$).) Then~\eqref{eq:32} can be re-expressed as:
\begin{equation}
  \label{eq:33}
  \begin{split}
    \log\, \norm{\dual{s}{s'}}^2
    &=\frac{\sqrt{-1}}{\pi}
    \int_X  f_1\,c_1(\rho)\,\log \norm{s'}
    +\frac{\sqrt{-1}}{\pi}
    \int_X f_2 \log\,\norm{s}\, c_1(\rho') \\
    & +\frac{\sqrt{-1}}{\pi} \int_X \log\,\norm{s}\:
    \d f_2 \wedge \dc \log\,\norm{s'} \\
    & +\frac{\sqrt{-1}}{\pi} \int_X \log\,\norm{s'} \:
    \d f_1 \wedge \dc \log\,\norm{s}\,,
  \end{split}
\end{equation}
a variation of \cite[6.5.1]{MR89b:32038}, due to O.~Gabber.
Observe that the integrations above make sense as ordinary smooth
differential forms, since, as a consequence of the assumptions on
$f_1$, $f_2$, both $\d f_1$ and $\d f_2$ vanish in a neighborhood
of $\supp E\cup \supp E'$.

\subsection{Pairing for good line bundles}
\label{sec:pairing-good-line}

Let us now include the divisor $D=p_1+\dots +p_N$ in the picture.
Accordingly, $\Bar{\sheaf{L}}$ and $\Bar{\sheaf{L}}'$ are good
hermitian line bundles on $X$ along $D$. Recall
$X_\varepsilon=X\setminus \bigcup_{i=1}^N \Delta_\varepsilon
(p_i)$.
The $2$-forms appearing in the integrand on the right hand side
of~\eqref{eq:33} are good forms along $D$, namely they belong
to $E^2_X(1)\lAngle D\rAngle_\good \eqdef \Gamma
(X,\she[2]{X}(1)\lAngle D\rAngle_\good)$.
\begin{definition}
  \label{def:8}
  With $X$, $D$, $\Bar{\sheaf{L}}$ and $\Bar{\sheaf{L}}'$ as
  above, set
  \begin{equation*}
    \log\, \norm{\dual{s}{s'}}^2
    = \lim_{\varepsilon\to 0} \int_{X_\varepsilon}
    \bigl(\text{Integrand $2$-form
        in RHS of~\eqref{eq:33}}\bigr)
  \end{equation*}
\end{definition}
The following is a direct generalization of \cite[Thm.
5.5]{math.CV/0211055} and \cite[Thm. 7.3]{bry:quillen}.
\begin{theorem}
  \label{thm:2}
  The cup-product in Hermitian-holomorphic Deligne cohomology
  with growth ($?=\good$ or $\pre$) extends the
  the norm on the determinant of cohomology line to good
  hermitian line bundles. Namely, let $\omega$ be a $2$-form
  corresponding to $[\Bar{\sheaf{L}}]\cup [\Bar{\sheaf{L}}']$ as
  in Lemma~\ref{lem:3}. Then
  \begin{equation}
    \label{eq:34}
    \lim_{\varepsilon\to 0}
    \frac{\sqrt{-1}}{\pi}\int_{X_\varepsilon} \omega
    = \log\, \norm{\dual{\mathbf{1}}{\mathbf{1}}}^2
  \end{equation}
  where the right hand side of~\eqref{eq:34} is interpreted in
  the sense of definition~\ref{def:8}, and $\mathbf{1}$ is the
  rational section of $\sheaf{L}$ and $\sheaf{L}'$ induced by the
  function $1$.
\end{theorem}
\begin{proof}
  We will proceed as in the proof of the corresponding result in
  ref.\ \cite{math.CV/0211055}.
  
  Let again $E$ and $E'$ the divisors of $\sheaf{L}$ and
  $\sheaf{L}'$. Recall their supports are disjoint. Hence we have
  an open cover $\cover{U}_X = \lbrace U_1,U_2\rbrace$, where
  $U_1 = X\setminus \supp E'$ and $U_2=X\setminus \supp E$.
  Moreover, choosing $f_1$, $f_2$ as in
  sect.~\ref{sec:determ-cohom} above, yields a partition of unity
  subordinated to the cover $\cover{U}_X$.
  
  Since $\cover{U}_X$ has two elements only, the only non-zero
  terms of the cocycle in the proof of Lemma~\ref{lem:3}
  representing $[\Bar{\sheaf{L}}]\cup [\Bar{\sheaf{L}}']$, are
  $\omega^0_i$ and $\omega^1_{ij}$, where the index $i$ takes the
  values $1,2$.  Furthermore, it is immediately verified that the
  $2$-form
  \begin{equation}
    \label{eq:35}
    \omega = f_1\,\omega^{0}_1 +f_2\,\omega^{1}_2
    +\d f_2 \wedge \omega^{1}_{21} 
  \end{equation}
  is a well-defined element of $E^2_X(1)\lAngle D\rAngle_?$ and
  represents the same cohomology class. Note that here $?=\good$
  \emph{and} $\pre$.
  
  Now, let $s$ and $s'$ be rational sections as
  in~\ref{sec:determ-cohom}.  With respect to the cover
  $\cover{U}_X$ they correspond to two pairs of rational
  functions, $s=\lbrace s_1,s_2\rbrace$, and $s'=\lbrace
  s'_1,s'_2\rbrace$. In fact $s'_1$ and $s_2$ are invertible in
  their respective domains of definition.  Let us now actually
  use that $s$ and $s'$ are the rational section $\mathbf{1}$, so
  that $s_2=1$ and $s'_1=1$.  With these choices we have:
  \begin{align*}
    \omega^0_1 &= \onehalf \log\rho_1 \d\dc \log \norm{s'}\\
    \omega^0_2 &= \log\norm{s} \: c_1(\rho')\\
    \omega^0_{21} &= \log \abs{s_1}\, \dc \log \norm{s'}
    -\dc \log\abs{s_1}\,\log\norm{s'}
  \end{align*}
  Then simple integrations by parts analogous to those in
  \cite[Thm. 5.5]{math.CV/0211055} yield:
  \begin{equation}
    \label{eq:36}
    \begin{split}
      \int_{X_\varepsilon} \omega & = \int_{X_\varepsilon}
      \bigl(\text{Integrand $2$-form
        in RHS of~\eqref{eq:33}}\bigr)\\
      &+ \onefourth \int_{\partial X_{\varepsilon}} f_1\, \bigl(
      \log\rho_1\, \dc \log\rho'_1
      -\log\rho'_1\, \dc \log\rho_1 \bigr)\,.
    \end{split}
  \end{equation}
  Since the metrics are good, it easily checked that each term in
  the second integral in~\eqref{eq:36} has an estimate of the
  type:
 \begin{equation*}
   \int_0^{2\pi} \abs{f_1(\varepsilon e^{\sqrt{-1}t})}
   \abs{\alpha(\varepsilon e^{\sqrt{-1}t})}
   \frac{\log\abs{\log{\varepsilon}}}{\abs{\log \varepsilon}}
   \, \d t
 \end{equation*}
 where $\alpha$ is bounded (and so is $f_1$). Hence, the second
 integral in~\eqref{eq:36} goes to zero as $\varepsilon\to 0$, as
 wanted.
\end{proof}

\section{Extremal hyperbolic metrics}
\label{sec:extr-hyperb-metr}

\subsection{Preliminaries}
\label{sec:preliminaries}

In this section we consider in more detail the setup of
example~\ref{ex:hyp-curves}: recall that $\d s^2_\hyp$ is a
hyperbolic metric on $U=X\setminus D$ satisfying the
condition~\eqref{eq:28}. The good extension of $T_U$ to $X$ is
$\sheaf{L}=T_X(-D)$. Thus we denote by $\Bar{\sheaf{L}}_\hyp$ the
pair $(\sheaf{L},\d s^2_\hyp)$.

More generally, we consider the space of good conformal metrics
$\d s^2$ on $X$. Namely $\d s^2$ is a conformal metric on $U$
locally written as $\d s^2= \rho \abs{\d z}^2$, where $\rho$
satisfies the same asymptotic conditions~\eqref{eq:28}. From
those conditions it follows that $\d s^2$ is good. We denote
$\cm{X}_\good$ the space of good conformal metrics on $X$, and
$\Bar{\sheaf{L}}$ will denote the pair $(\sheaf{L},\d s^2)$.

Note that $\cm{X}_\good$ is an affine space over
$\Gamma(X,\she[0]{X})$.  Namely, if $\d s^2$ and $\d s'^2\in
\cm{X}$, there exists a smooth function $\sigma :X\to \RR$ such
that $\d s'^2 = e^{\sigma}\d s^2$. The hyperbolic metric $\d
s^2_\hyp$ belongs to $\cm{X}_\good$.

Let $V$ be an adapted neighborhood of $p\in D$, that is, a disk
$\Delta_a$ of radius $a>0$ centered at $p$, where we write
\begin{math}
  \d s^2 = \rho \abs{\d z}^2.
\end{math}
Using the non vanishing section $s=z\del/\del z$ of
$\sheaf{L}\vert_V$ let us also set:
\begin{equation}
  \label{eq:38}
  \Tilde\rho (z) \eqdef \norm{s}^2 \coin
  \d s^2 \Bigl(z\frac{\del}{\del z},
  z\frac{\del}{\del z}\Bigr)\,,
\end{equation}
so that $\rho = \Tilde\rho / \abs{z}^2$ and
$\Tilde\rho$ has the asymptotic behavior:
\begin{equation}
  \label{eq:37}
  \log\Tilde\rho (z) = -\log
    \bigl( \log^2 \abs{z} \bigr)+\alpha (z)
\end{equation}
at $z=0$.  If, on the other hand, $(V,z)$ is a neighborhood of an
ordinary point $p\in U$, then we simply set $\Tilde\rho=\rho$ and
$s=\del/\del z$.

Finally, let $\omega_\rho$ be the K\"ahler form associated to $\d
s^2$. Locally we have
\begin{equation*}
  \omega_\rho = \ihalf \rho \d z \wedge \d\Bar z\,,
\end{equation*}
and the associated \emph{area} is simply the integral
\begin{math}
  A_\rho (X) = \int_X \omega_\rho\,.
\end{math}

\subsection{The Liouville equation}
\label{sec:liouville-equation}

The hyperbolic metric $\d s^2_\hyp$ whose scalar curvature is
equal to $-1$ is of course of particular importance from the
point of view of uniformization theory.  With respect to a local
coordinate $z$ the constant negative curvature condition
translates into the nonlinear PDE
\begin{equation}
  \label{eq:39}
  \frac{\del^2\log\rho_\hyp}{\del z\del\Bar z} =
  \onehalf \rho_\hyp\,,
\end{equation}
known as the \emph{Liouville equation,} which must be
supplemented by the asymptotic condition~\eqref{eq:28} at the
points of $D$. The equation can be recast in the more invariant
fashion:
\begin{equation}
  \label{eq:40}
  c_1(\rho_\hyp) = \sqrt{-1}\,\omega_{\rho_\hyp}\,.
\end{equation}
Note that $c_1(\rho_\hyp)=c_1(\Tilde\rho_\hyp)$, and the latter
computes $c_1(\sheaf{L})$.

It is well-known that~\eqref{eq:39} is the extremum of an
energy-type functional. A global construction for
eq.~\eqref{eq:40} was given in~\cite{zogtak1987-1} for curves of
type $(0,n)$ (with $n\geq 3$) and
in~\cite{zogtak1987-2,math.CV/0204318} for curves of type $(g,0)$
(with $g\geq 2$), and more generally, in the latter case,
in~\cite{math.CV/0211055}.

The following is a direct generalization of~\cite[Thm. 5.1 and
Cor.  5.6]{math.CV/0211055}:
\begin{theorem}
  \label{thm:3}
  Denote by $\log\,
  \norm{\dual{\Bar{\sheaf{L}}}{\Bar{\sheaf{L}}}}$ the integral of
  the good $2$-form determined by $[\Bar{\sheaf{L}}]\cuphat
  [{\Bar{\sheaf{L}}}]$, as in Theorem~\ref{thm:2} and consider
  the quantity:
  \begin{equation}
    \label{eq:41}
    S (\d s^2) \eqdef 
    \log\, \norm{\dual{\Bar{\sheaf{L}}}{\Bar{\sheaf{L}}}}
    + \frac{1}{2\pi} A_\rho(X)\,.
  \end{equation}
  The metric $\d s^2$ satisfies equation~\eqref{eq:40} (that is:
  $\d s^2 =\d s^2_\hyp$) if and only if:
  \begin{equation*}
    \frac{\d}{\d t}\Bigr\rvert_{t=0}S (e^{t\sigma}\d s^2)=0\,, 
  \end{equation*}
  for any $\sigma\in\Gamma(X,\she[0]{X})$.
\end{theorem}
\begin{proof}
  First, fix an open cover $\cover{U}_X=\lbrace V_i\rbrace_{i\in
    I}$ with adapted neighborhoods. Let $z_i$ be the
  corresponding coordinate on $V_i$. (We will drop the \cech
  index $i$ if there is no danger of confusion.) Assume the same
  conditions on the cover and the integration map explained in
  sect.~\ref{sec:trace} before Lemma~\ref{lem:5}.
  \begin{lemma}
    \label{lem:4}
    If $\sheaf{L}$ and $\d s^2$ are trivialized by sections $s_i$
    on $V_i$ as in sect.~\ref{sec:preliminaries}, the resulting
    cocycle belongs to
    \begin{math}
      \Tot^2 \vC{\cover{U}_X}{\ndhh{X,?}{1}},
    \end{math}
    with $?=\good$ or $\pre$.
  \end{lemma}
  \begin{proof}
    This follows at once from the definitions.
  \end{proof}
  Using $?=\pre$, the cup square in pre-log-log Hermitian
  holomorphic Deligne cohomology has the form~\eqref{eq:23} (with
  $\rho'=\rho$), and let again $\Omega$ be the cocycle
  representing the resulting image of the class
  $[\Bar{\sheaf{L}}]\cuphat [\Bar{\sheaf{L}}]$ in
  $\HHH^2(X,\she{X}(1)\lAngle D\rAngle_?)$ as per
  Lemma~\ref{lem:3}. Using sect.~\ref{sec:trace} and
  Lemma~\ref{lem:5}, eq.~\eqref{eq:41} is written as:
  \begin{equation*}
    S(\d s^2) = 
    \frac{\sqrt{-1}}{2\pi}\dual{\Omega}{\Sigma_\mathbf{T}}
    + \frac{1}{2\pi} A_\rho (X)
  \end{equation*}
  In order to carry out the calculation of
  \begin{math}
    \frac{\d}{\d t}\Bigr\rvert_{t=0}S (e^{t\sigma}\d s^2),
  \end{math}
  we can use the arguments of
  \cite{math.CV/0211055,math.CV/0204318} provided~\eqref{eq:44}
  is taken into account when integrating by parts. Since the
  calculation is the same we will not reproduce it here. As a
  result, we have:
  \begin{equation}
    \label{eq:45}
    \begin{split}
      \frac{\d}{\d t}\Bigr\rvert_{t=0}S (e^{t\sigma}\d s^2)
      &= \frac{\sqrt{-1}}{2\pi}
      \int_X \sigma
      \Bigl( c_1(\sheaf{L}) -\sqrt{-1}\omega_\rho \Bigr) \\
      &- \frac{\sqrt{-1}}{8\pi}
      \lim_{\varepsilon\to 0}
      \sum_{p_i\in D}
      \int_{C_\varepsilon(p_i)}
      \bigl(\log\Tilde\rho_i \dc\sigma
      -\sigma \dc\log \Tilde\rho_i 
      \bigr)\,,
  \end{split}
\end{equation}
where $C_\varepsilon(p_i)$ is a circle of radius $\varepsilon$
around the point $p_i$. The first integral in eq.~\eqref{eq:45}
must be understood in the sense of eq.~\eqref{eq:42}. The second
integral is treated as the analogous one at the end of the proof
of Theorem~\ref{thm:2}. Then its limit as $\varepsilon\to 0$ is
zero. This finishes the proof.
\end{proof}
\begin{remark}
  The use of explicit correction terms, as in ref.\
  \cite{zogtak1987-1}, is not needed in the present framework.
  Notice that on an adapted neighborhood $V$ the section
  $\del/\del z$ ought to be considered as a \emph{rational}
  section of $\sheaf{L}\vert_V$: we have $\del/\del z = z^{-1}s$.
  In the language of ref.\ \cite{math.AG/0404122} the
  corresponding forms will be \emph{mixed,} that is they present
  both \emph{log-log} and \emph{log}-type singularities, thereby
  losing local integrability.
\end{remark}
The following lemma provides a precise comparison. Let $p\in
\supp D$ and let us identify the adapted neighborhood $(V,z)$
with a disk $\Delta_a$ of radius $a$ centered at $z=0$. Let
$A(\varepsilon,a)=\Delta_a\cap X_\varepsilon$ be the resulting
annulus of radii $a$ and $\varepsilon$. Recall that $\rho (z) =
\norm{\del/\del z}^2$ and $\Tilde\rho (z) = \norm{s}^2$. Then we
have:
\begin{lemma}
  \label{lem:6}
  \begin{equation*}
    \int_{A(\varepsilon,a)} \onehalf \log\Tilde\rho \,
    \delb\del \log\Tilde\rho
    = \int_{A(\varepsilon,a)} \onehalf
    \del\log\rho \wedge \delb\log\rho
    - \tate \log \varepsilon
    - \tate \log (\log^2 \varepsilon) + O(1)\,.
  \end{equation*}
\end{lemma}
\begin{proof}
  One has the identity
  \begin{equation*}
    \log\Tilde\rho \, \delb\del \log\Tilde\rho
    = \del\log\rho \wedge \delb\log\rho
    + \onehalf \d \Bigl\lbrace
    \log\Tilde\rho\, \dc \log\Tilde\rho
    +\log \abs{z}^2 \dc\log\abs{z}^2
    -2 \log\Tilde\rho\, \dc \log\abs{z}^2
    \Bigr\rbrace\,.
  \end{equation*}
  The estimate on the boundary term follows from a direct
  calculation and the fact that $\log\Tilde\rho\, \dc
  \log\Tilde\rho$ is a good form.
\end{proof}

\bibliography{general.bib} \bibliographystyle{hamsplain}

\end{document}